\documentclass[12pt]{article}
\usepackage{amsmath}
\usepackage{graphicx}
\pdfoutput=1
\usepackage{graphics}
\usepackage{CJK}
\usepackage{epsfig}
\usepackage{multirow}
\usepackage{lipsum}
\usepackage{makecell}

\usepackage{subfigure}
\usepackage{enumerate}
\usepackage{diagbox}
\usepackage{caption}
\usepackage{float}
\input{amssym.tex}

\def\bc{\begin{center}}       \def\ec{\end{center}}
\def\ba{\begin{array}}        \def\ea{\end{array}}
\def\be{\begin{equation}}     \def\ee{\end{equation}}
\def\bea{\begin{eqnarray}}    \def\eea{\end{eqnarray}}
\def\beaa{\begin{eqnarray*}}  \def\eeaa{\end{eqnarray*}}

\def\mathbb{\Bbb}

\begin{document}

\baselineskip 18pt
\centerline {\bf \large Bifurcation of limit cycles for a class of cubic Hamiltonian
}
\vskip 0.1 true cm
\centerline {\bf \large   systems with nesting period annuli }

\vskip 0.3 true cm

\centerline{\bf  Yuan Chang, Liqin Zhao$^{*}$, Qiuyi Wang}
 \centerline{ School of Mathematical Sciences, Beijing Normal University,} \centerline{Beijing 100875, The People's Republic of China}

\footnotetext[1]{This work was supported by National Natural Science Foundation of China (12071037) *Corresponding author.
E-mail:  zhaoliqin@bnu.edu.cn (L. Zhao); yuanchang1996@163.com   }
\vskip 0.2 true cm

\noindent{\bf Abstract}~~In this paper, we obtain the upper bound of the number of zeros of Abelian integral $I(h)=\oint_{\Gamma_{h}}g(x,y)dy-f(x,y)dx$, where $\Gamma_{h}$ is the closed oribit defined by
$$H(x,y)=x^2-y^2+ax^4+cy^4+bx^2y^2,(a,b,c)\in\mathbb{R}^3, c\neq0,$$
and $f(x,y)$, $g(x,y)$ are polynomials in $(x,y)$ with degree $n$. Furthermore, we consider the Hopf and homoclinic bifurcation when $a=-1$, $b=-2$, $c=1$ and $n=3$, and obtain 18 distributions in which system has at least $3$ limit cycles for each case.

\section{Introduction and the main results}

Many periodic phenomena in nature can be described by limit cycle, which has a lot of applications in economics \cite{puu}, physics \cite{3body}, astronomy \cite{collin} and so on. As one of the main topics of qualitative theory, the maximum number of limit cycles for planar continuous differential systems is widely studied by researchers. For more results about this topic, we prefer the reader to see [3, 5--8, 12, 14, 19--22] and reference therein. In this paper, we deal with the limit cycles of a class of cubic Hamiltonian systems with the origin being a saddle under polynomial perturbations.

Consider the near-Hamiltonian system
$$\dot{x}=H_{y}(x,y)+\varepsilon f(x, y), \quad  \dot{y}=-H_{x}(x,y)+\varepsilon g(x, y),\eqno(1.1)$$
where $0<|\varepsilon|\ll1$, $H(x,y)$ is a polynomial of degree $m+1$, $f(x,y)$ and $g(x,y)$ are polynomials of degree $n$. Abelian integral of system (1.1) is $$I(h)=\oint_{\Gamma_{h}}g(x,y)dx-f(x,y)dy,$$
where $\Gamma_{h}$ is a family of closed orbits defined by $H(x,y)=h$, $h\in\Sigma$, $\Sigma$ is an union of several open intervals on which $\Gamma_{h}$ is defined. The question of finding the upper bound $Z(m,n)$ of the number of isolated zeros of $I(h)$ is called the weakened Hilbert's 16th problem \cite{Arnold}. By Poincar\'{e}-Pontryagin Theorem, the number of zeros of $I(h)$ provides an upper bound for the number of limit cycles of system (1.1) if $I(h)$ is not identically zero.

For the general results of finding $Z(m,n)$, Khovansky \cite{khovan} and Varchenko \cite{varchen} proved independently the finiteness of $Z(m,n)$, but they did not give its explicit expression. For $m=2$, the normal form of Hamiltonian
$$H(x,y)=\frac{1}{2}(x^2+y^2)-\frac{1}{3}x^3+axy^2+\frac{1}{3}by^3$$
has been given if the unperturbed system has annulus, Horozov and Iliev \cite{horo} derived $Z(2,n)\leq5n+15$ by analysing the corresponding Picard--Fuchs equations. For $m\geq3$, Petrov [19,20] estimated the number of zeros of $I(h)$ for Hamiltonians
$$H(x,y)=y^2+x^2-x^4,\quad H(x,y)=y^2-2x^2+x^4.$$
Zhao and Zhang \cite{zhaoyulin1999} studied the following Hamiltonian
$$H(x,y)=\frac12y^2+U(x), \deg U(x)=4,\eqno(1.2)$$
they gave four kinds of normal forms for (1.2) having at least one center
and derived the upper bound $7n+5$ for the number of isolated zeros of $I(h)$. Liu \cite{cjliu8} estimated the least upper bound of the number of zeros of $I(h)$ for
$$H(x,y)=y^2+(x^4-x^2-\lambda x), -\frac{2\sqrt{2}}9<\lambda<0.$$
Zhou and Li [30,31] gave the upper bound of the number of isolated zeros of $I(h)$ for
$$H(x,y)=-x^2+x^4+y^4$$
and
$$H(x,y)=x^2+y^2+ax^4+bx^2y^2+cy^4, ac(b^2-4ac)\neq0$$
with $a>0$, $b=0$ and $c=1$.
Wu, Zhang and Li \cite{Wujuanjuan} proved the number of isolated zeros of $I(h)$ does not exceed $2\left[\frac{n-1}{4}\right]+12\left[\frac{n-3}{4}\right]+23$ for
$$H(x,y)=x^2+y^2-x^4+ax^2y^2+y^4,a>2.$$

Recently, Yang and Zhao \cite{yangjihua2017} proved the number of isolated zeros of $I(h)$ does not exceed $48n+12$ for Hamiltonian
$$H(x,y)=-x^2+ax^4+bx^2y^2+cy^4, c\neq0$$
except the butterfly phase portrait, and $51n+622$ for each period annulus of butterfly phase portrait in \cite{yangbf}. Yang \cite{yang2019ds} considered the Hamiltonian
$$H(x,y)=ax^2+by^2+cx^2y^2-x^4+y^4,$$
and proved the number of isolated zeros of $I(h)$ does not exceed $90n+24$. Chen and Yu \cite{chenyu} discussed the following Hamiltonian
$$H(x,y)=x^2+y^2+ax^4+bx^2y^2+cy^4, ac(b^2-4ac)\neq0,$$
and proved $I(h)$ has at most $45n+323$ zeros for each period annulus, taking into account the multiplicity.
Motivated by [2, 26, 27, 28] , we consider a class of cubic Hamiltonian with the following form
$$H(x,y)=x^2-y^2+ax^4+cy^4+bx^2y^2,(a,b,c)\in\mathbb{R}^3, c\neq0,\eqno(1.3)$$
and the corresponding vector field is
$$
\begin{cases}
\dot{x}=2y(-1-bx^2+2cy^2),\\
\dot{y}=-2x(1+2ax^2+by^2).
\end{cases}\eqno(1.4)
$$
The phase portrait of system (1.4) are summarized as listed in Fig. 1 and Fig. 2.

Consider the following perturbed system
$$
\begin{cases}
\dot{x}=2y(-1-bx^2+2cy^2)+\varepsilon f(x,y),\\
\dot{y}=-2x(1+2ax^2+by^2)+\varepsilon g(x,y),
\end{cases}\eqno(1.5)
$$
where $0<|\varepsilon|\ll1$, $f(x,y)$ and $g(x,y)$ are polynomials in $x$ and $y$ with degree $n$. Our main results are the following two theorems.
\begin{figure}[htbp]\centering
\centering
\includegraphics[width=12cm,height=9.8cm]{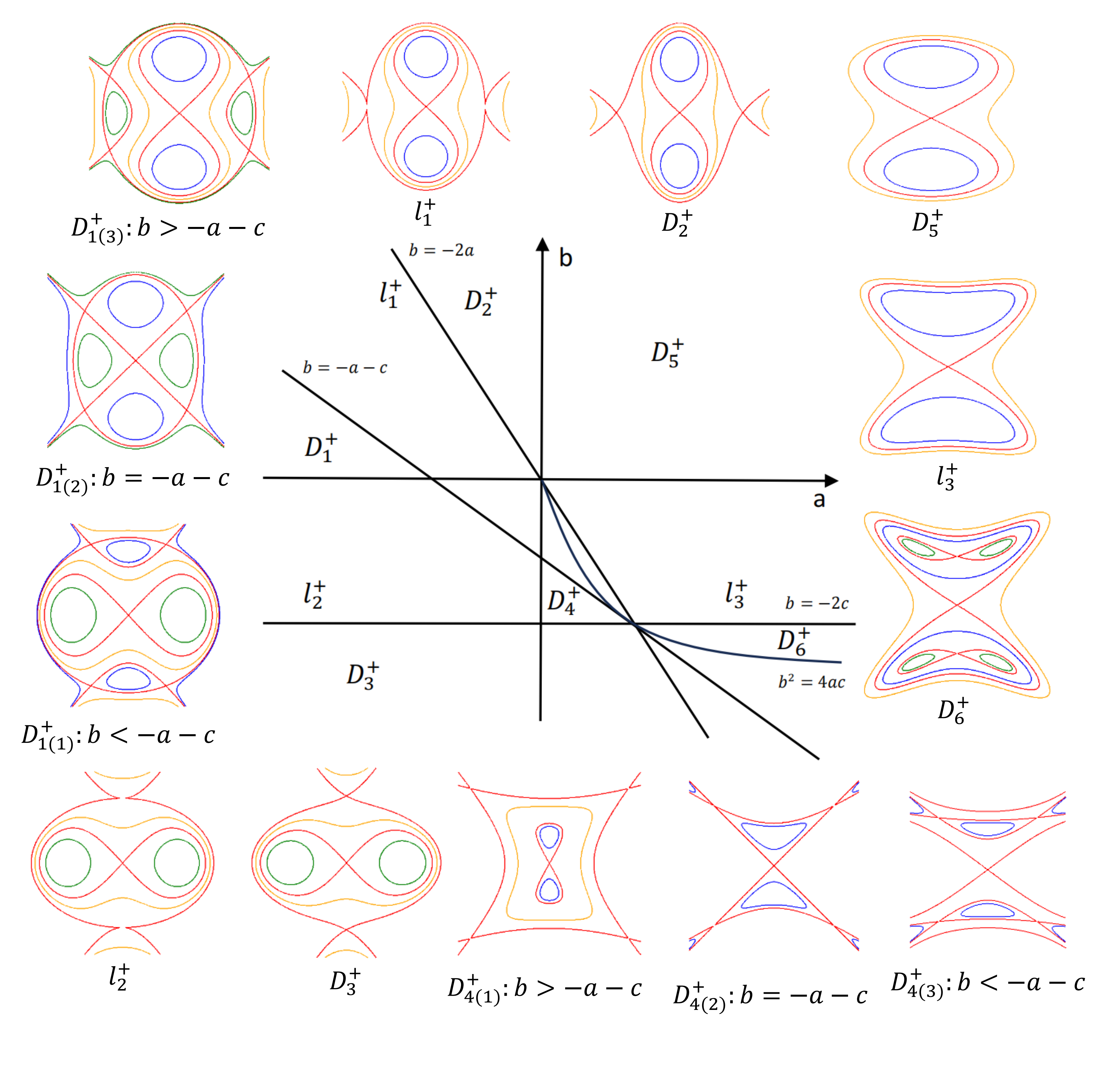}
\caption{\small Phase portraits of system (1.4) for $c>0$. $D_{1}^{+}=\{a<0,b+2a<0,b+2c<0\}$, $D_{2}^{+}=\{a<0,b+2a>0\}$, $l_{1}^{+}=\{a<0,b+2a=0\}$, $D_{3}^{+}=\{a<0,b+2c<0\}$, $l_{2}^{+}=\{a<0,b+2c=0\}$, $D_{4}^{+}=\{a\geq0,b<0,b+2c>0,b^2-4ac>0\}$,  $D_{5}^{+}=\{a\geq0,b\geq0\}\cup\{b^2-4ac\leq0,b+2c>0\}$, $l_{3}^{+}=\{b^2-4ac<0,b+2c=0\}$, $D_{6}^{+}=\{b^2-4ac<0,b+2c<0\}$. System has no period annuli for $(a,b)\in\mathbb{R}^2\setminus\{D_{1}^{+}\cup D_{2}^{+}\cup D_{3}^{+}\cup D_{4}^{+}\cup D_{5}^{+}\cup D_{6}^{+}\cup l_{1}^{+}\cup l_{2}^{+}\cup l_{3}^{+}\}$.}
\end{figure}
\begin{figure}[htbp]\centering
\centering
\includegraphics[width=12cm,height=8.5cm]{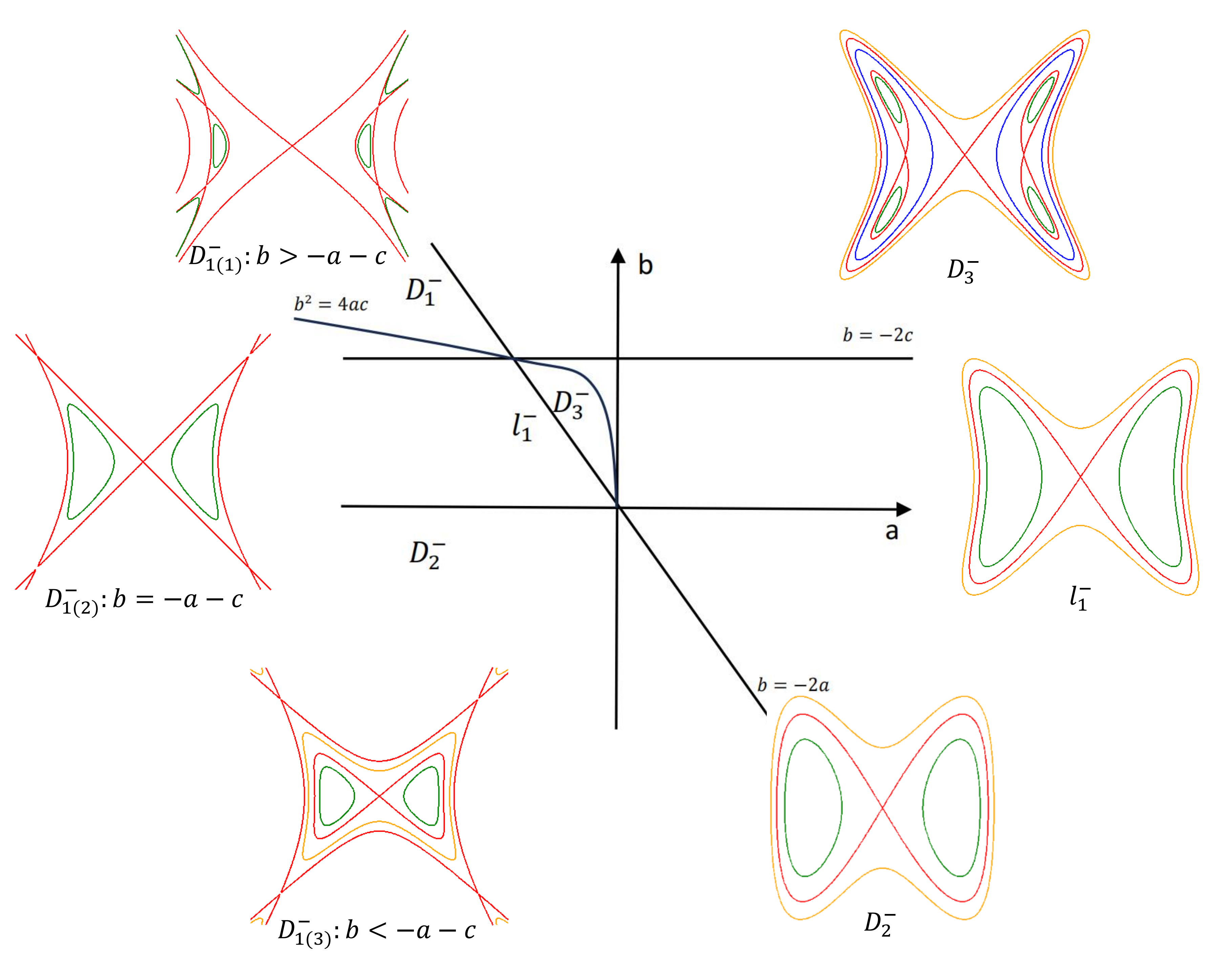}
\caption{\small Phase portraits of system (1.4) for $c<0$. $D_{1}^{-}=\{b^2-4ac>0,b+2c>0,b+2a<0\}$, $D_{2}^{-}=\{a<0,b<0\}\cup\{b^2-4ac\leq0,b+2a<0\}$, $l_{1}^{-}=\{b^2-4ac<0,b+2a=0\}$, $D_{3}^{-}=\{b^2-4ac<0,b+2a>0\}$. System has no period annuli for $(a,b)\in\mathbb{R}^2\setminus\{D_{1}^{-}\cup D_{2}^{-}\cup D_{3}^{-}\cup l_{1}^{-}\}$.}
\end{figure}
\begin{table}[htbp]
\centering
\resizebox{\linewidth}{!}{
\begin{tabular}{ccc|cc}
 \hline
  \makebox[0.2\textwidth][c]{\multirow{2}{*}{Case}}&\makebox[0.2\textwidth][c]{\multirow{2}{*}{Former result}}&\makebox[0.2\textwidth][c]{\multirow{2}{*}{Our result}}&\makebox[0.2\textwidth][c]{\multirow{2}{*}{Case}}&\makebox[0.2\textwidth][c]{\multirow{2}{*}{Our result}} \\[10pt]
   \hline
\multirow{2}{*}{$D^{+}_{1(1)}$}    &\multirow{2}{*}{$90n+24$}&\multirow{2}{*}{$58n+121$} &\multirow{2}{*}{$D^{+}_{4(3)}(a=0)$} &\multirow{2}{*}{$45n-58$} \\[10pt]
 \multirow{2}{*}{$D^{+}_{1(2)}$} &\multirow{2}{*}{$86n+13$}&\multirow{2}{*}{$54n + 109$}  &\multirow{2}{*}{$D^{+}_{5}(a\neq0)$} &\multirow{2}{*}{$31n+69$} \\[10pt]
 \multirow{2}{*}{$D^{+}_{1(3)}$}    &\multirow{2}{*}{$90n+24$} &\multirow{2}{*}{$58n+121$} &\multirow{2}{*}{$D^{+}_{5}(a=0)$} &\multirow{2}{*}{$49n-60$} \\[10pt]
 \multirow{2}{*}{$l^{+}_{1}$} &\multirow{2}{*}{$47n+16$}&\multirow{2}{*}{$31n + 66$}  &\multirow{2}{*}{$l^{+}_{3}$} &\multirow{2}{*}{$31n+68$} \\[10pt]
  \multirow{2}{*}{$D^{+}_{2}$}    &\multirow{2}{*}{$47n+16$}&\multirow{2}{*}{$31n + 67$}  &\multirow{2}{*}{$D^{+}_{6}$} &\multirow{2}{*}{$208n + 1089$} \\[10pt]
   \multirow{2}{*}{$l^{+}_{2}$}    &\multirow{2}{*}{$47n+16$} &\multirow{2}{*}{$31n + 66$}  & \multirow{2}{*}{$D^{-}_{1(1)}$} &\multirow{2}{*}{$27n+55$} \\[10pt]
    \multirow{2}{*}{$D^{+}_{3}$}    &\multirow{2}{*}{$47n+16$} &\multirow{2}{*}{$31n + 67$}  &\multirow{2}{*}{$D^{-}_{1(2)}$} &\multirow{2}{*}{$27n+55$} \\[10pt]
    \multirow{2}{*}{$D^{+}_{4(1)}(a\neq0)$}    &\multirow{2}{*}{-} &\multirow{2}{*}{$31n+66$}  &\multirow{2}{*}{$D^{-}_{1(3)}$} &\multirow{2}{*}{$31n+66$} \\[10pt]
    \multirow{2}{*}{$D^{+}_{4(2)}(a\neq0)$}    &\multirow{2}{*}{-} &\multirow{2}{*}{$27n+55$}  &\multirow{2}{*}{$D^{-}_{2}$} &\multirow{2}{*}{$31n+69$} \\[10pt]
    \multirow{2}{*}{$D^{+}_{4(3)}(a\neq0)$}    &\multirow{2}{*}{-} &\multirow{2}{*}{$27n+55$}  &\multirow{2}{*}{$l^{-}_{1}$} &\multirow{2}{*}{$31n+68$} \\[10pt]
    \multirow{2}{*}{$D^{+}_{4(1)}(a=0)$}    &\multirow{2}{*}{-} &\multirow{2}{*}{$49n-60$}  &\multirow{2}{*}{$D^{-}_{3}$} &\multirow{2}{*}{$208n + 1089$} \\[10pt]
    \multirow{2}{*}{$D^{+}_{4(2)}(a=0)$}    &\multirow{2}{*}{-} &\multirow{2}{*}{$45n-58$}  &\ &\  \\[10pt]
   \hline
 \end{tabular}}
\multirow{4}{*}{Table 1. The upper bound of the number of zeros of $I(h)$ for each case in Fig 1 and Fig 2}
 \end{table}

\vskip 0.2 true cm
{ \it \noindent{\bf Theorem 1.1.}~Suppose that $(a,b,c)\in\mathbb{R}^3$ with $c(b^2-4ac)\neq0$, then for arbitrary polynomial $f(x,y)$ and $g(x,y)$ with degree $n$, we obtain the number of isolated zeros of $I(h)$ of system (1.5) does not exceed $208n+1089$ (taking into account multiplicity).}

Then for system (1.5) with $a=-1$, $b=-2$, $c=1$ and $n=3$, we have
$$
\begin{cases}
\dot{x}=2y(-1+2x^2+2y^2)+\varepsilon \sum\limits_{i+j=0}^{3}a_{ij}x^iy^j,\\
\dot{y}=-2x(1-2x^2-2y^2)+\varepsilon \sum\limits_{i+j=0}^{3}b_{ij}x^iy^j.
\end{cases}\eqno(1.6)
$$
System (1.6)$|_{\varepsilon=0}$ has two centers $(\pm\frac{1}{\sqrt{2}},0)$, and a double-homoclinic loop through saddle $(0,0)$, see $l_{2}^{+}$ in Fig 1. Consider the Hopf and homoclinic bifurcations of system (1.6), we obtain

\vskip 0.2 true cm
{ \it \noindent{\bf Theorem 1.2.}~For Hopf and homoclinic bifurcations, we obtain that system (1.6) has

\noindent(i)~at most $3$ limit cycles near each center, and has

\noindent(ii)~at most $5$ limit cycles near the double-homoclinic loop.

\noindent Furthermore, we obtain $18$ distributions in which system (1.6) has at least $3$ limit cycles for each case.
}

{\bf\noindent Remark 1.}~(i)~By Poincar\'{e}-Pontryagin Theorem, the number of zeros of $I(h)$ provides an upper bound for the number of limit cycles of system (1.5).

\noindent(ii)~The upper bound of the number of zeros of $I(h)$ for all the cases in Fig 1 and Fig 2 are summarized as listed in Table 1. By decreasing the degree of three polynomials in second-order differential operator, we obtain a lower upper bound for $D_{1}$, $D_{2}$, $D_{3}$, $l_{1}^{+}$ and $l_{2}^{+}$.

This paper is organized as follows. In section 2, we derive the algebraic structure of $I(h)$ for system (1.5), which is generated by nine generators $I_{01},\ I_{03},\ I_{21},\ I_{23},\ I_{12}$, $I_{11}$, $I_{13}$, $I_{02}$ and $I_{22}$, where
$$I_{ij}=\oint_{\Gamma_{h}}x^{i}y^{j}dx.$$
In section 3, we derive four Picard-Fuchs equations for $\{I_{01}, I_{03}, I_{21}, I_{23}\}$,  $\{I_{11}, I_{13}\}$ and $\{I_{02}, I_{22}\}$, $\{I_{01}, I_{03}, I_{21}, I_{23}, I_{12}\}$ respectively, and then we show that each of
$$\omega_{1}=\frac{Z^{\prime}}{I^{\prime}_{01}},\quad \omega_{2}=\frac{I_{13}}{I_{11}},\quad \omega_{3}=\frac{I_{22}}{I_{02}},\quad \bar{\omega}_{1}=\frac{\bar{Z}}{I_{01}}$$
satisfies a Riccati equation, where
$$Z=-\frac{b+2c}{6(b^2-4ac)}I_{03}+\frac{b+2a}{2(b^2-4ac)}I_{21}+\frac{1}{3}I_{23}, \quad \bar{Z}=\frac{2}{5}I_{03}+\frac{6c}{5b}I_{21}.$$
Using these results, we shell prove there exists three second-order differential operator $L_{i}(h), i=1,2$ such that $L_{1}(h)I(h)=F_{1}(h)$, $L_{2}(h)L_{1}(h)I(h)=F_{2}(h)$, and can estimate the number of zeros of $I(h)$ from the number of zeros of $F_{i}(h),i=1,2$, where $F_{i}(h)$ is analytic. The results in theorem 1.1 are proved by the ideas from [9, 15, 29, 32]. In section 4, we consider the Hopf and homoclinic bifurcation of $l_{2}^{+}$ with asymptotic expansions of Abelian Integrals by the ideas from [10,17], and discuss the coexistence of limit cycles.

Throughout the paper, we denote by $\#\{f(h)=0,h\in(s,t)\}$ the number of isolated zeros of $f(h)$ on $(s,t)$ taking into account the multiplicity, denote by $\#\{I(h)=0, (a,b,c)\in D\}$ the number of isolated zeros of $I(h)$ for case $D$ in Fig. 1 and Fig. 2, taking into account the multiplicity, and denote by $A^T$ the transpose of matrix $A$. We also use the following notations:
$$
\begin{aligned}
V_{1}(h)&=(I_{01}(h),I_{03}(h),I_{21}(h),I_{23}(h))^T,\ V_{2}(h)=(I_{11}(h),I_{13}(h))^T, \\
V_{3}(h)&=(I_{02}(h),I_{22}(h))^T,\  V_{4}(h)=(I_{01}(h),I_{03}(h),I_{21}(h),I_{23}(h),I_{12}(h))^T,\\
\Theta_{1}(h)&=(f_{1}(h), f_{2}(h), f_{3}(h), f_{4}(h))V_{1}(h),\ \Theta_{2}(h)=(g_{1}(h),g_{2}(h))V_{2}(h), \\
\Theta_{3}(h)&=(l_{1}(h),l_{2}(h))V_{3}(h), \Theta_{4}(h)=(f_{1}(h), f_{2}(h), f_{3}(h), f_{4}(h), f_{5}(h))V_{4}(h).\\
\end{aligned}
$$

\section{The algebraic structure of $I(h)$}

\vskip 0.2 true cm
{\it \noindent{\bf Lemma 2.1.}~For $ac(b^2-4ac)\neq0$, if $i+j>4 (i+j=n)$ is odd, then
$$
I_{ij}(h)=\bar{f}_{1}(h)I_{01}+\bar{f}_{2}(h)I_{03}+\bar{f}_{3}(h)I_{21}+\bar{f}_{4}(h)I_{23}+\bar{f}_{5}(h)I_{12}, \eqno(2.1)$$\\
and
$$
I_{2i+1,2j+1}(h)=\bar{g}_{1}(h)I_{11}+\bar{g}_{2}(h)I_{13}, \quad if\  i+j>1 (i+j=\frac{n}{2}-1),\eqno(2.2)$$
$$
I_{2i,2j}(h)=\bar{l}_{1}(h)I_{02}+\bar{l}_{2}(h)I_{22}, \quad if\  i+j>1 (i+j=\frac{n}{2}),\eqno(2.3)
$$
where the coefficients of $I_{ij}$ in (2.1)--(2.3) are real polynomials in $h$ with
$$
\begin{aligned}
{\rm deg} \bar{f}_{1}(h)&\leq\left[\frac{n-1}{4}\right], {\rm deg} \bar{f}_{2}(h), \bar{f}_{3}(h),\bar{f}_{5}(h)\leq\left[\frac{n-3}{4}\right], {\rm deg} \bar{f}_{4}(h)\leq\left[\frac{n-1}{4}\right]-1, \\
&{\rm deg} \bar{g}_{1}(h),\bar{l}_{1}(h)\leq\left[\frac{n-2}{4}\right], {\rm deg} \bar{g}_{2}(h), \bar{l}_{2}(h)\leq\left[\frac{n}{4}\right]-1.
\end{aligned}
$$
}
\vskip 0.2 true cm
\noindent{\bf Proof}.~Without loss of generality, we only prove (2.1) here, (2.2) and (2.3) can be proved similarly,
Multiplying (1.3) by $x^{i}y^{j-4}dx$ and integrating on $\Gamma_{h}$ give

$$I_{ij}=\frac{1}{c}(hI_{i,j-4}+I_{i,j-2}-I_{i+2,j-4}-aI_{i+4,j-4}-bI_{i+2,j-2}).\eqno(2.4)$$
Derivating (1.3) with respect to $x$, we obtain
$$x-y\frac{\partial y}{\partial x}+2ax^3+bxy^2+bx^2y\frac{\partial y}{\partial x}+2cy^3\frac{\partial y}{\partial x}=0.\eqno(2.5)$$
Multiplying (2.5) by $x^{i-3}y^{j}dx$ and integrating on $\Gamma_{h}$ give
$$I_{ij}=\frac{1}{2a}(-\frac{i-3}{j+2}I_{i-4,j+2}+\frac{2c(i-3)}{j+4}I_{i-4,j+4}-I_{i-2,j}+\frac{(i-j-3)b}{j+2}I_{i-2,j+2}).\eqno(2.6)$$
(2.4) and (2.6) can be reduced to
$$I_{ij}=\frac{j}{2c(i+j+1)}(2hI_{i,j-4}-I_{i+2,j-4}+\frac{i+2j-3}{j-2}I_{i,j-2}-\frac{(i+j+1)b}{j-2}I_{i+2,j-2}),\eqno(2.7)$$
$$I_{ij}=\frac{1}{2a(i+j+1)}(2(i-3)hI_{i-4,j}-(2i+j-2)I_{i-2,j}+\frac{j(i-3)}{j+2}I_{i-4,j+2}-\frac{(i+j+1)jb}{j+2}I_{i-2,j+2}).\eqno(2.8)$$
When $n=5$, (2.7) and (2.8) give
$$I_{05}=\frac{5}{12c}(2hI_{01}-\frac{7}{3}I_{03}-I_{21}-2bI_{23}),\qquad I_{14}=\frac{1}{3c}(3I_{12}-3bI_{32}),$$
$$I_{41}=\frac{1}{12a}(2hI_{01}+\frac{1}{3}I_{03}-7I_{21}-2bI_{23}),\qquad I_{32}=\frac{1}{12a}(-6I_{12}-3bI_{14}),$$
which yield the result for $n=5$. Suppose that the result hold for $n=2k-1(k\geq3)$, which implies
$$I_{ij}=\tilde{f}_{1}(h)I_{01}+\tilde{f}_{2}(h)I_{03}+\tilde{f}_{3}(h)I_{21}+\tilde{f}_{4}(h)I_{23}+\tilde{f}_{5}(h)I_{12},$$
where deg $\tilde{f}_{1}(h)\leq\left[\frac{2k-2}{4}\right]$, deg $\tilde{f}_{2}(h), \tilde{f}_{3}(h), \tilde{f}_{5}(h)\leq\left[\frac{2k-4}{4}\right]$, deg $\tilde{f}_{4}(h)\leq\left[\frac{2k-2}{4}\right]-1$.
Then for $n=2k+1$, taking $(i,j)=(0,2k+1), (1,2k), (2,2k-1)$  in (2.7) and $(i,j)=(3,2k-2), (4,2k-3), \ldots, (2k,1)$ in (2.8), respectively, we obtain
$$
A\left(
\begin{matrix}
I_{2k,1}\\
I_{2k-1,2}\\
\vdots\\
I_{4,2k-3}\\
I_{3,2k-1}\\
I_{2,2k-1}\\
I_{1,2k}\\
I_{0,2k+1}\\
\end{matrix}
\right)=\frac{1}{4k+4}\left(\begin{matrix}
\frac{1}{a}[2(2k-3)hI_{2k-4,1}-(4k-1)I_{2k-2,1}+\frac{2k-3}{3}I_{2k-4,3}]\\
\frac{1}{a}[2(2k-4)hI_{2k-5,2}-(4k-2)I_{2k-3,2}+(k-2)I_{2k-5,4}]\\
\vdots\\
\frac{1}{a}[2hI_{0,2k-3}-(2k+3)I_{2,2k-3}+\frac{2k-3}{2k-1}I_{0,2k-1}]\\
\frac{1}{a}[-(2k+2)I_{1,2k-2}]\\
\frac{1}{c}[2(2k-1)hI_{2,2k-5}-(2k-1)I_{4,2k-5}+\frac{(4k-3)(2k-1)}{2k-3}I_{2,2k-3}]\\
\frac{1}{c}[4khI_{1,2k-4}-2kI_{3,2k+4}+\frac{(4k-2)2k}{2k-2}I_{1,2k-2}]\\
\frac{1}{c}[2(2k+1)hI_{0,2k-3}-(2k+1)I_{2,2k-3}+\frac{(4k-1)(2k+1)}{2k-1}I_{0,2k-1}]\\
\end{matrix}\right)
,$$
where
$$A=\left(
\begin{matrix}
1&0&\frac{b}{6a}&0&\cdots&0&0&0&0&0\\
0&1&0&\frac{b}{4a}&\cdots&0&0&0&0&0\\
\vdots&\vdots&\vdots&\vdots&\vdots&\vdots&\vdots&\vdots&\vdots&\vdots\\
0&0&0&0&\cdots&1&0&\frac{(2k-3)}{2(2k-1)a}&0&0\\
0&0&0&0&\cdots&0&1&0&\frac{(2k-2)b}{4ak}&0\\
0&0&0&0&\cdots&\frac{(2k-1)b}{2(2k-3)c}&0&1&0&0\\
0&0&0&0&\cdots&0&\frac{2kb}{2(2k-2)c}&0&1&0\\
0&0&0&0&\cdots&1&0&\frac{(2k+1)b}{2(2k-1)c}&0&1\\
\end{matrix}\right).
$$
Since det$A=\frac{(4ac-b^2)^2}{16a^2c^2}\neq0$, thus
\begin{align*}
I_{ij}=&h(\hat{f}_{1}(h)I_{01}+\hat{f}_{2}(h)I_{03}+\hat{f}_{3}(h)I_{21}+\hat{f}_{4}(h)I_{23}+\hat{f}_{5}(h)I_{12})+\tilde{f}_{1}(h)I_{01}+\tilde{f}_{2}(h)I_{03}\\
&+\tilde{f}_{3}(h)I_{21}+\tilde{f}_{4}(h)I_{23}+\tilde{f}_{5}(h)I_{12}\\
=&\bar{f}_{1}(h)I_{01}+\bar{f}_{2}(h)I_{03}+\bar{f}_{3}(h)I_{21}+\bar{f}_{4}(h)I_{23}+\bar{f}_{5}(h)I_{12},
\end{align*}
where
\begin{align*}
{\rm deg} \hat{f}_{1}(h)\leq\left[\frac{2k-4}{4}\right], &{\rm deg} \hat{f}_{2}(h), \hat{f}_{3}(h), \hat{f}_{5}(h)\leq\left[\frac{2k-6}{4}\right], {\rm deg} \hat{f}_{4}(h)\leq\left[\frac{2k-4}{4}\right]-1,\\
{\rm deg} \bar{f}_{1}(h)\leq\left[\frac{2k}{4}\right], &{\rm deg} \bar{f}_{2}(h), \bar{f}_{3}(h), \bar{f}_{5}(h)\leq\left[\frac{2k-2}{4}\right], {\rm deg} \bar{f}_{4}(h)\leq\left[\frac{2k}{4}\right]-1.
\end{align*}
This end the proof. $\diamondsuit$

By lemma 2.1, we can get the following lemma.

{\it \noindent{\bf Lemma 2.2.}~For $ac(b^2-4ac)\neq0$, $I(h)$ can be expressed as
$$
\begin{aligned}
I(h)=&[f_{1}(h)I_{01}+f_{2}(h)I_{03}+f_{3}(h)I_{21}+f_{4}(h)I_{23}+f_{5}(h)I_{12}]\\
&+[g_{1}(h)I_{11}+g_{2}(h)I_{13}]+[l_{1}(h)I_{02}+l_{2}(h)I_{22}],
\end{aligned}
$$
where
\begin{align*}
{\rm deg} f_{1}(h)&\leq\left[\frac{n-1}{4}\right], {\rm deg} f_{2}(h), f_{3}(h), f_{5}(h)\leq\left[\frac{n-3}{4}\right], {\rm deg} f_{4}(h)\leq\left[\frac{n-1}{4}\right]-1, \\
&{\rm deg} g_{1}(h),l_{1}(h)\leq\left[\frac{n-2}{4}\right], {\rm deg} g_{2}(h),l_{2}(h)\leq\left[\frac{n}{4}\right]-1.
\end{align*}
}
\section{Proof of Theorem 1.1}

In this section, we will mainly deal with the case $D_{6}^{+}$ in Fig. 1. The system (1.4) has seven singular points:

\noindent$\bullet$ Four elementary centers $C_{1,2,3,4}(\pm\sqrt{\frac{b+2c}{b^2-4ac}},\pm\sqrt{-\frac{b+2a}{b^2-4ac}})$,\\[10pt]

\vspace{-1.5em}
\noindent$\bullet$ Three elementary saddles $S_{0}(0,0)$, $S_{1,2}(0,\pm\sqrt{\frac{1}{2c}})$.

The period annulus only around $C_{i}$ (i=1,\ 2,\ 3,\ 4) corresponds to $h\in\Sigma_{1}$, the period annulus expanding from the double homoclinic orbits through saddle points $S_{0}$ to the double homoclinic orbits through saddle points $S_{1}$ and $S_{2}$ corresponds $h\in\Sigma_{2}$, the period annulus expand outwards from the double homoclinic orbits through saddle points $S_{0}$ corresponds to $h\in\Sigma_{3}$,  where
$$\Sigma_{1}=(\frac{a+b+c}{b^2-4ac}, -\frac{1}{4c}), \Sigma_{2}=(-\frac{1}{4c}, 0), \Sigma_{3}=(0, +\infty).\eqno(3.1)$$

From the discussion in section 2, we know that
$$
\begin{aligned}
I(h)=&[f_{1}(h)I_{01}+f_{2}(h)I_{03}+f_{3}(h)I_{21}+f_{4}(h)I_{23}+f_{5}(h)I_{12}]\\
&+[g_{1}(h)I_{11}+g_{2}(h)I_{13}]+[l_{1}(h)I_{02}+l_{2}(h)I_{22}],\ {\rm for}\ h\in\Sigma_{1},
\end{aligned}
$$
$$
I(h)=[f_{1}(h)I_{01}+f_{2}(h)I_{03}+f_{3}(h)I_{21}+f_{4}(h)I_{23}]+[l_{1}(h)I_{02}+l_{2}(h)I_{22}],\ {\rm for}\ h\in\Sigma_{2},
$$
and
$$
I(h)=f_{1}(h)I_{01}+f_{2}(h)I_{03}+f_{3}(h)I_{21}+f_{4}(h)I_{23},\ {\rm for}\ h\in\Sigma_{3}.
$$
In order to estimate the number of zeros of $I(h)$ for $h\in\Sigma_1\cup\Sigma_2\cup\Sigma_3$, we first introduce the following useful result in \cite{Gavrilov} and \cite{zhao2021}. Let $V$ be a finite-dimensional vector space of functions, real-analytic on an open interval $\mathbb{I}$.

\vskip 0.2 true cm
{\it \noindent{\bf Definition 3.1.}~(\cite{Gavrilov}).
We say that $V$ is a Chebyshev space, provided that each non-zero function in $V$ has at most dim$(V)-1$ zeros, counted with multiplicity.}

Let $S$ be the solution space of a second order linear analytic differential equation
$$x''+a_1(t)x'+a_2(t)x=0 \eqno(3.2)$$
on an open interval $\mathbb{I}$.

\vskip 0.2 true cm
{\it \noindent{\bf Proposition 3.1.}~(\cite{Gavrilov}).
The solution space S of (3.2) is a Chebyshev space on the interval $\mathbb{I}$ if
and only if there exists a nowhere vanishing solution $x_{0}(t)\in S\ (x_{0}(t)\neq0,\forall t\in \mathbb{I})$.}

\vskip 0.2 true cm
{\it \noindent{\bf Proposition 3.2.}~(\cite{Gavrilov}).
Suppose the solution space of the homogeneous equation (3.2) is a
Chebyshev space and let $R(t)$ be an analytic function on $\mathbb{I}$ having $l$ zeros (counted with multiplicity). Then every solution $x(t)$ of the non-homogeneous equation
$$x''+a_1(t)x'+a_2(t)x=R(t)$$
has at most $l+2$ zeros on $\mathbb{I}$.}

From proposition 3.1 and 3.2, Zhou et al.\cite{zhao2021} gave the following proposition.

\vskip 0.2 true cm
{\it \noindent{\bf Proposition 3.3.}~(lemma 6 of \cite{zhao2021}). Assume that $f(t)=f_{1}(t)+f_{2}(t)$, where $f$, $f_{1}$, $f_{2}$ are real--analytic functions in $t\in(\alpha,\beta)$. If there exist polynomials $p_{1}(t)$, $p_{2}(t)$, $p_{3}(t)$ respectively such that
$$L(t)f(t)=r(t), \quad L(t)f_{1}(t)=0,$$
where $L(t)=p_{3}(t)\frac{d^2}{dt^2}+p_{2}(t)\frac{d}{dt}+p_{1}(t)$, then we have
$$
\begin{aligned}
\#\{f(t)=0,t\in(\alpha,\beta)\}\leq&3\#\{p_{3}(t)=0, t\in(\alpha,\beta)\}+3\#\{f_{1}(t)=0, t\in(\alpha,\beta)\}\\
&+\#\{r(t)=0,t\in(\alpha,\beta)\}+2.
\end{aligned}
$$}

\vskip 0.2 true cm
{\it \noindent{\bf Lemma 3.1.}
For $ac(b^2-4ac)\neq0$, the vector functions $V_{i}(h),\ i=1,2,3,4$ satisfying the following Picard-Fuchs equations
$$V_{i}=(A_{i}h+B_{i})V_{i}^{\prime},\ i=1,2,3,4,
\eqno(3.3)$$
where $$A_1=\left( \begin {array}{cccc} 2&0&0&0\\ \noalign{\medskip}\frac{3}{4c}
&1&0&0\\ \noalign{\medskip}-\frac{1}{4a}&0&1&0\\ \noalign{\medskip}a_{11}&a_{12}&a_{13}&\frac23\end {array}
 \right),\qquad
B_1=\left( \begin {array}{cccc} 0&\frac13&-1&0\\ \noalign{\medskip}0&\frac{3}{8c}&-\frac{3}{8c}&-\frac{b}{4c}-\frac12\\ \noalign{\medskip}0&-\frac{1}{24a}&\frac{3}{8a}&\frac{b}{12a}+\frac16\\ \noalign{\medskip}0&b_{12}&b_{13}&b_{14}\end {array}
 \right),$$
$$
A_2=\left( \begin{array}{cc} \frac43&0\\ \noalign{\medskip}{-\frac {4(b+2a)}{5(b^2-4ac)}}
&\frac45\end{array} \right),
\qquad
B_2=\left( \begin{array}{cc} \frac{1}{3a}&\frac{b}{9a}+\frac29\\ \noalign{\medskip}-{\frac {b+2a}{5
a(b^2-4ac)}}&-\frac{12ac+b^2+16ab+16a^2}{15a(b^2-4ac)}\end{array} \right),
$$
$$
A_3=\left( \begin{array}{cc} \frac43&0\\ \noalign{\medskip}{\frac {4(b+2c)}{15(b^2-4ac)}}
&\frac45\end{array} \right),
\qquad
B_3=\left( \begin{array}{cc} \frac{1}{3c}&-\frac{b}{3c}-\frac23\\ \noalign{\medskip}{\frac {b+2c}{15
c(b^2-4ac)}}&-\frac{b^2+16bc+16c^2+12ac}{15c(b^2-4ac)}\end{array} \right),
$$
$$
A_4=\left( \begin{array}{cc} A_{1}&0\\ \noalign{\medskip}0
&1\end{array} \right),
\qquad
B_4=\left( \begin{array}{cc} B_{1}&0\\ \noalign{\medskip}0&-\frac{a+b+c}{b^2-4ac}\end{array} \right)
$$}
with
$$\begin{aligned}
a_{11}=&\frac{4ac+ab+bc}{8ac(b^2-4ac)}, \quad\qquad\qquad\qquad\qquad a_{12}=\frac{b+2c}{6(b^2-4ac)},\\
a_{13}=&-\frac{b+2a}{2(b^2-4ac)}, \qquad\qquad\qquad\qquad b_{12}=\frac{8ac+3ab+bc}{48ac(b^2-4ac)},\\
b_{13}=&-\frac{8ac+ab+3bc}{16ac(b^2-4ac)}, \qquad b_{14}=-\frac{24abc+20a^2c+ab^2+20ac^2+b^2c}{24ac(b^2-4ac)}.
\end{aligned}$$

\vskip 0.2 true cm
\noindent{\bf Proof}.~From (1.3), we have
$$\frac{\partial y}{\partial h}=\frac{1}{2y(2cy^2+bx^2-1)},$$
which implies
$$I_{ij}^{\prime}=\frac{j}{2}\oint_{\Gamma_{h}}\frac{x^{i}y^{j-2}}{2cy^2+bx^2-1}dx.\eqno(3.4)$$
Hence,
$$I_{ij}=\oint_{\Gamma_{h}}\frac{x^{i}y^{j}(2cy^2+bx^2-1)}{2cy^2+bx^2-1}dx=\frac{2b}{j+2}I_{i+2,j+2}^{\prime}+\frac{4c}{j+4}I_{i,j+4}^{\prime}-\frac{2}{j+2}I^{\prime}_{i,j+2},\eqno(3.5)$$
Multiplying both side of (3.4) by h, we obtain
$$
\begin{aligned}
hI_{ij}^{\prime}&=\frac{j}{2}\oint_{\Gamma_{h}}\frac{x^{i}y^{j-2}(x^2-y^2+ax^4+bx^2y^2+cy^4)}{2cy^2+bx^2-1}dx\\
&=I_{i+2,j}^{\prime}-\frac{j}{j+2}I_{i,j+2}^{\prime}+aI_{i+4,j}^{\prime}+\frac{jb}{j+2}I_{i+2,j+2}^{\prime}+\frac{jc}{j+4}I_{i,j+4}^{\prime}.
\end{aligned}
\eqno(3.6)$$
By Green formula, we can get
$$
\begin{aligned}
I_{ij}&=\oint_{\Gamma_{h}}x^iy^jdx=-\frac{j}{i+1}\oint_{\Gamma_{h}}x^{i+1}y^{j-1}dy\\
&=-\frac{j}{i+1}\oint_{\Gamma_{h}}x^{i+1}y^{j-1}\frac{-2x(1+2ax^2+by^2)}{2y(-1+bx^2+2cy^2)}dx\\
&=\frac{2}{i+1}I_{i+2,j}^{\prime}+\frac{4a}{i+1}I_{i+4,j}^{\prime}+\frac{2bj}{(i+1)(j+2)}I_{i+2,j+2}^{\prime}.
\end{aligned}\eqno(3.7)
$$
From (3.5)--(3.7), we have
$$I_{ij}=\frac{2}{i+j+1}(2hI_{ij}^{\prime}-I_{i+2,j}^{\prime}+\frac{j}{j+2}I_{i,j+2}^{\prime}).\eqno(3.8)$$
By (3.8), we have the following equations

$$\begin{aligned}I_{01}&=2hI^{\prime}_{01}+I^{\prime}_{21}-\frac13I^{\prime}_{03},\qquad I_{03}=hI^{\prime}_{03}-\frac12I^{\prime}_{23}+\frac3{10}I^{\prime}_{05},\\
I_{21}&=hI^{\prime}_{21}-\frac12I^{\prime}_{41}+\frac1{6}I^{\prime}_{23},\qquad I_{23}=\frac23hI^{\prime}_{23}-\frac13I^{\prime}_{43}+\frac15I^{\prime}_{25},\\
I_{12}&=hI_{12}^{\prime}-\frac12I_{32}^{\prime}+\frac14I_{14}^{\prime},\qquad I_{11}=\frac43hI^{\prime}_{11}+\frac13I^{\prime}_{31}-\frac19I^{\prime}_{13}.\\
I_{13}&=\frac45hI^{\prime}_{13}+\frac15I^{\prime}_{33}-\frac3{25}I^{\prime}_{15},\quad I_{02}=\frac{4}{3}hI_{02}^{\prime}+\frac{1}{3}I_{22}^{\prime}-\frac{1}{6}I_{04}^{\prime}\\
I_{22}&=\frac{4}{5}hI_{22}^{\prime}+\frac{1}{5}I_{42}^{\prime}-\frac{1}{10}I_{24}^{\prime}.\end{aligned}$$
Then we can get the conclusion by (2.7) and (2.8). This ends the proof. $\diamondsuit$

\vskip 0.2 true cm
{\it \noindent{\bf Lemma 3.2.}
For $ac(b^2-4ac)\neq0$, let
$$Z(h)=-\frac{b+2c}{6(b^2-4ac)}I_{03}(h)+\frac{b+2a}{2(b^2-4ac)}I_{21}(h)+\frac{1}{3}I_{23}(h),$$
then $I_{01}, I_{03}, I_{21}, I_{12}, Z$ satisfy the equation
$$
G_{1}(h) \left(
\begin{matrix}
I_{01}^{\prime\prime}\\
I_{03}^{\prime\prime}\\
I_{21}^{\prime\prime}\\
Z^{\prime\prime}\\
I_{12}^{\prime\prime}\\
\end{matrix}
\right)=\left(
\begin{matrix}
d_{11}(h) \quad d_{12}(h)\\
d_{21}(h) \quad d_{22}(h)\\
d_{31}(h) \quad d_{32}(h)\\
d_{41}(h) \quad d_{42}(h)\\
0 \quad\quad\quad 0\\
\end{matrix}
\right)\left(
\begin{matrix}
I_{01}^{\prime}\\
Z^{\prime}\\
\end{matrix}
\right),\eqno(3.9)
$$
where
$$\begin{aligned}
G_{1}(h)&=\frac{1}{12ac}h(4ah+1)(4ch+1)(h-\frac{a+b+c}{b^2-4ac}),\\
d_{11}(h)&=-\frac{1}{12ac}h(8ach+a+c)(h-\frac{a+b+c}{b^2-4ac}),\\
d_{12}(h)&=-\frac{1}{24ac}(2abh+8ach+2bch+a+b+c),\\
d_{21}(h)&=-\frac{1}{8ac}h(4ah+1)(h-\frac{a+b+c}{b^2-4ac}),\\
\end{aligned}$$
$$
\begin{aligned}
d_{22}(h)&=\frac{1}{8ac}h(4ah+1)(b+2c),\\
d_{31}(h)&=\frac{1}{24ac}h(4ch+1)(h-\frac{a+b+c}{b^2-4ac}),\\
d_{32}(h)&=-\frac{1}{24ac}h(4ch+1)(b+2a),\\
d_{41}(h)&=-\frac{h}{24ac(b^2-4ac)}(2abh+8ach+2bch+a+b+c)(h-\frac{a+b+c}{b^2-4ac}),\\
d_{42}(h)&=\frac{1}{12ac}h(8ach+a+c)(h-\frac{a+b+c}{b^2-4ac}).
\end{aligned}$$}
\noindent{\bf Proof}.~Derivating (3.3) with respect to h, we get
$$(A_{i}h+B_{i})V_{i}^{\prime\prime}=(I-B_{i})V_{i},$$
where $I$ is a $4\times4$ identity matrix. By elementary manipulations, we can get the conclusion. This ends the proof. $\diamondsuit$

\noindent{\bf Remark 2.}~Noting that (3.9) can not be used directly for $h\in\tilde{\Sigma}$ if $G(h)$ has a zero in it, which divides $\tilde{\Sigma}$ into $\tilde{\Sigma}_{1}$ and $\tilde{\Sigma}_{2}$. Then the number of zeros of $f(h)$ on $\tilde{\Sigma}$ does not exceed the number of zeros of $f(h)$ on $\tilde{\Sigma}_{1}\cup\tilde{\Sigma}_{2}$ plus 1, where $f(h)$ is a polynomial in $h$. This issue can be found in
$$S=\{D_{2}^{+}\cup D_{3}^{+}\cup D_{5}^{+}\cup l_{3}^{+}\cup D_{6}^{+}\cup D_{3}^{-}\cup l_{1}^{-}\cup D_{2}^{-}\}.$$ Hence, we can get the upper bound of the number of zeros of $I(h)$ for all period annuli of each case in $S$.

By lemma 3.1 and elementary manipulations, we can get the following lemma.

\vskip 0.2 true cm
{\it \noindent{\bf Lemma 3.3.}
Set $\omega_{1}=\frac{Z^{\prime}}{I_{01}^{\prime}}$, $\omega_{2}=\frac{I_{13}}{I_{11}}$, $\omega_{3}=\frac{I_{22}}{I_{02}}$, then we obtain

\noindent(i)~$$G_{1}(h)\omega_{1}^{\prime}=-d_{12}(h)\omega_{1}^{2}+(d_{42}(h)-d_{11}(h))\omega_{1}+d_{41}(h),$$
where $G_{1}(h)$, $d_{11}(h)$, $d_{12}(h)$, $d_{41}(h)$, $d_{42}(h)$ are defined in lemma 3.2.

\noindent(ii)~$$G_{2}(h)\omega_{2}^{\prime}=-a_{2}(h)\omega_{2}^{2}+(a_{4}(h)-a_{1}(h))\omega_{2}+a_{3}(h),$$
where $G_{2}(h)=\frac{4}{15a}(4ah+1)(h-\frac{a+b+c}{b^2-4ac})$, and
$$\begin{aligned}
a_{1}(h)&=\frac45h-\frac{12ac+b^2+16ab+16a^2}{15a(b^2-4ac)},\quad a_{2}(h)=-\frac{b}{9a}-\frac29,\\
a_{3}(h)&=\frac{4(b+2a)}{5(b^2-4ac)}h+\frac{b+2a}{5a(b^2-4ac)},\quad a_{4}(h)=\frac43h+\frac{1}{3a}.
\end{aligned}$$

\noindent(iii)~$$G_{3}(h)\omega_{3}^{\prime}=-b_{2}(h)\omega_{3}^{2}+(b_{4}(h)-b_{1}(h))\omega_{3}+b_{3}(h),$$
where $G_{3}(h)=\frac{4}{15c}(4ch+1)(h-\frac{a+b+c}{b^2-4ac})$, and
$$\begin{aligned}
b_{1}(h)&=\frac45h-\frac{12ac+b^2+16bc+16c^2}{15c(b^2-4ac)},\quad b_{2}(h)=\frac{b}{3c}-\frac23,\\
b_{3}(h)&=-\frac{4(b+2c)}{15(b^2-4ac)}h-\frac{b+2c}{15c(b^2-4ac)},\quad b_{4}(h)=\frac43h+\frac{1}{3c}.
\end{aligned}$$
}

\noindent{\bf Proof}.~By lemma 3.1 and lemma 3.2, we can get the conclusion. This ends the proof. $\diamondsuit$

\vskip 0.2 true cm
{\it \noindent{\bf Lemma 3.4.}For $h\in\Sigma_{1}$, there exist polynomials $R_{3}(h)$, $R_{2}(h)$ and $R_{1}(h)$ of $h$ with degree respectively $n_{1}$, $n_{1}-1$ and $n_{1}-2$ such that $L_{1}(h)\Theta_{3}(h)=0$, where $n_{1}=\left[\frac{n-2}{4}\right]+\left[\frac{n}{4}\right]+2$, and
$$L_{1}(h)=R_{3}(h)\frac{d^{2}}{dh^{2}}+R_{2}(h)\frac{d}{dh}+R_{1}(h).$$

}

\noindent{\bf Proof.}~We first assert that
$$
\begin{cases}
\Theta_{3}(h)=l_{11}(h)I_{02}^{\prime\prime}(h)+l_{12}(h)I_{22}^{\prime\prime}(h),\\
\Theta_{3}^{\prime}(h)=l_{21}(h)I_{02}^{\prime\prime}(h)+l_{22}(h)I_{22}^{\prime\prime}(h),\\
\Theta_{3}^{\prime\prime}(h)=l_{31}(h)I_{02}^{\prime\prime}(h)+l_{32}(h)I_{22}^{\prime\prime}(h),
\end{cases}
$$
where
$$
\begin{aligned}
{\rm deg} l_{11}&(h)\leq\left[\frac{n-2}{4}\right]+2, {\rm deg} l_{12}(h)\leq\left[\frac{n}{4}\right]+1, {\rm deg} l_{21}(h)\leq\left[\frac{n-2}{4}\right]+1,\\
&{\rm deg} l_{22}(h)\leq\left[\frac{n}{4}\right], {\rm deg} l_{31}(h)\leq\left[\frac{n-2}{4}\right], {\rm deg} l_{32}(h)\leq\left[\frac{n}{4}\right]-1.
\end{aligned}
$$
In fact, it follows from lemma 3.1 that
$$
\begin{aligned}
\Theta_{3}(h)&=(l_{1}(h), l_{2}(h))V_{3}=(l_{1}(h), l_{2}(h))(A_{3}h+B_{3})(I-A_{3})^{-1}(A_{3}h+B_{3})V_{3}^{\prime\prime}\\
&=l_{11}(h)I_{11}^{\prime\prime}(h)+l_{12}(h)I_{13}^{\prime\prime}(h),
\end{aligned}
$$
where $I$ is a $2\times2$ identity matrix.
So we can calculate deg$l_{11}(h)\leq\left[\frac{n-2}{4}\right]+2$ and deg$l_{12}(h)\leq\left[\frac{n}{4}\right]+1$ from lemma 3.1. For $\Theta_{3}^{\prime}(h)$, we have
$$
\begin{aligned}
\Theta_{3}^{\prime}(h)&=(l_{1}(h), l_{2}(h))V^{\prime}_{3}(h)+(l^{\prime}_{1}(h), l^{\prime}_{2}(h))V_{3}(h)\\
&=[(l_{1}(h),l_{2}(h))+(l_{1}^{\prime}(h),l_{2}^{\prime}(h))(A_{3}h+B_{3})](I-A_{3})^{-1}(A_{3}h+B_{3})V_{3}^{\prime\prime}(h).
\end{aligned}
$$
The result for $\Theta_{3}^{\prime\prime}(h)$ can be proved similarly.

Next, suppose that
$$R_{1}(h)=\sum_{i=0}^{n_{1}-2}r_{1i}h^i, R_{2}(h)=\sum_{j=0}^{n_{1}-1}r_{2j}h^j, R_{3}(h)=\sum_{k=0}^{n_{1}}r_{3k}h^k\eqno(3.10)$$
are polynomials of $h$ with coefficients $r_{1i}$, $r_{2j}$ and $r_{3k}$ to be determined such that $L_{1}(h)\Theta_{3}(h)=0$ for
$$0\leq i\leq n_{1}-2,\ 0\leq j\leq n_{1}-1,\ 0\leq k\leq n_{1}.\eqno(3.11)$$
By calculation, we have
$$
\begin{aligned}
L_{1}(h)\Theta_{3}(h)=&R_{3}(h)\Theta_{3}^{\prime\prime}(h)+R_{2}(h)\Theta_{3}^{\prime}(h)+R_{1}(h)\Theta_{3}(h)\\
=&R_{3}(h)[l_{31}(h)I_{02}^{\prime\prime}(h)+l_{32}(h)I_{22}^{\prime\prime}(h)]\\
&+R_{2}(h)[l_{21}(h)I_{02}^{\prime\prime}(h)+l_{22}(h)I_{22}^{\prime\prime}(h)]\\
&+R_{1}(h)[l_{11}(h)I_{02}^{\prime\prime}(h)+l_{12}(h)I_{22}^{\prime\prime}(h)]\\
=&[R_{3}(h)l_{31}(h)+R_{2}(h)l_{21}(h)+R_{1}(h)l_{11}(h)]I_{02}^{\prime\prime}(h)\\
&+[R_{3}(h)l_{32}(h)+R_{2}(h)l_{22}(h)+R_{1}(h)l_{12}(h)]I_{22}^{\prime\prime}(h)\\
=&P(h)I_{02}^{\prime\prime}(h)+Q(h)I_{22}^{\prime\prime}(h),
\end{aligned}
$$
where $P(h)$ and $Q(h)$ are polynomials of $h$ with degree no more than $n_{1}+\left[\frac{n-2}{4}\right]$ and $n_{1}+\left[\frac{n}{4}\right]-1$ respectively. Let
$$P(h)=\sum^{n_{1}+\left[\frac{n-2}{4}\right]}_{l=0}p_{l}h^{l},\quad Q(h)=\sum^{n_{1}+\left[\frac{n}{4}\right]-1}_{l=0}q_{l}h^{l},$$
where $p_{l}$ and $q_{l}$ are expressed by $r_{1i}$, $r_{2j}$ and $r_{3k}$ of (3.10) linearly, $i$, $j$ and $k$ satisfy (3.11). So $L_{1}(h)\Theta_{3}(h)=0$ is satisfied if we let
$$p_{l}=0,\quad q_{l}=0,\quad 0\leq i\leq n_{1}+\left[\frac{n-2}{4}\right],\quad 0\leq j\leq n_{1}+\left[\frac{n}{4}\right]-1.\eqno(3.12)$$
System (3.12) is a homogeneous linear equations with $2n_{1}+\left[\frac{n-2}{4}\right]+\left[\frac{n}{4}\right]+1$ equations about $3n_{1}$ variables of $r_{1i}$, $r_{2j}$ and $r_{3k}$ for $i$, $j$ and $k$ satisfy (3.11). It is easy to check if $n_{1}\geq\left[\frac{n-2}{4}\right]+\left[\frac{n}{4}\right]+2$ that
$$3n_{1}-\left(2n_{1}+\left[\frac{n-2}{4}\right]+\left[\frac{n}{4}\right]+1\right)\geq1.$$
Then we set $n_{1}=\left[\frac{n-2}{4}\right]+\left[\frac{n}{4}\right]+2$, which follows that from the theory of linear algebra that there exist $r_{1i}$, $r_{2j}$ and $r_{3k}$ such that (3.12) holds. Thus we get the conclusion of the existence of $L_{1}(h)$. This ends the proof. $\diamondsuit$

\vskip 0.2 true cm
{\it \noindent{\bf Lemma 3.5.}~For $h\in\Sigma_{1}$, we obtain $L_{1}(h)I(h)=F_{1}(h)$, where $L_{1}(h)$ is defined in lemma 3.4, and
$$\begin{aligned}F_{1}(h)=&\frac{1}{G^{2}_{2}(h)}[Q_{1}(h)I_{11}+Q_{2}(h)I_{13}]+\frac{1}{G_{1}(h)}[Q_{3}(h)I^{\prime}_{01}+Q_{4}(h)Z^{\prime}]\\
&+Q_{5}(h)I^{\prime}_{03}+Q_{6}(h)I^{\prime}_{21}+Q_{7}(h)I^{\prime}_{12}\end{aligned}$$
with
$$G_{2}(h)=\frac{4}{15a}(4ah+1)(h-\frac{a+b+c}{b^2-4ac}),$$
$$
\begin{aligned}
{\rm deg} Q_{1}(h)\leq& n_{1}+\left[\frac{n-2}{4}\right]+2, {\rm deg} Q_{2}(h)\leq n_{1}+\left[\frac{n}{4}\right]+1,\\
{\rm deg} Q_{3}(h)\leq& n_{1}+\left[\frac{n-1}{4}\right]+3, {\rm deg} Q_{4}(h)\leq n_{1}+\left[\frac{n-1}{4}\right]+2,\\
{\rm deg}& Q_{5}(h),Q_{6}(h),Q_{7}(h)\leq n_{1}+\left[\frac{n-3}{4}\right]-1.
\end{aligned}
$$
}
\noindent{\bf Proof.}~We first assert that

$$
\begin{cases}
\Theta_{2}(h)=g_{11}(h)I_{11}(h)+g_{12}(h)I_{13}(h),
\ \\
\Theta_{2}^{\prime}(h)=\frac{1}{G_{2}(h)}\left[g_{21}(h)I_{11}(h)+g_{22}(h)I_{13}(h)\right],
\ \\
\Theta_{2}^{\prime\prime}(h)=\frac{1}{G^{2}_{2}(h)}\left[g_{31}(h)I_{11}(h)+g_{32}(h)I_{13}(h)\right],
\end{cases}
\eqno(3.13)$$
where
$$
\begin{aligned}
{\rm deg} &g_{11}(h)\leq\left[\frac{n-2}{4}\right], {\rm deg} g_{12}(h)\leq\left[\frac{n}{4}\right]-1, {\rm deg} g_{21}(h)\leq\left[\frac{n-2}{4}\right]+1,\\
&{\rm deg} g_{22}(h)\leq\left[\frac{n}{4}\right], {\rm deg} g_{31}(h)\leq\left[\frac{n-2}{4}\right]+2, {\rm deg} g_{32}(h)\leq\left[\frac{n}{4}\right]+1.
\end{aligned}
$$
In fact, by lemma 3.1, we obtain
$$
\begin{aligned}
\Theta_{2}^{\prime}(h)=&(g_{1}(h),g_{2}(h))V^{\prime}_{2}(h)+(g_{1}^{\prime}(h),g_{2}^{\prime}(h))V_{2}(h)\\
=&[(g_{1}(h),g_{2}(h))(A_{2}h+B_{2})^{-1}+(g_{1}^{\prime}(h),g_{2}^{\prime}(h))]V_{2}(h)\\
=&\frac{1}{G_{2}(h)}\left[g_{21}(h)I_{11}(h)+g_{22}(h)I_{13}(h)\right],
\end{aligned}
$$
and
$$
\begin{aligned}
\Theta_{2}^{\prime\prime}(h)=&(g_{1}^{\prime\prime}(h),g_{2}^{\prime\prime}(h))V_{2}(h)+2(g_{1}^{\prime}(h),g_{2}^{\prime}(h))V_{2}^{\prime}(h)+(g_{1}(h),g_{2}(h))V_{2}^{\prime\prime}(h)\\
=&((g_{1}^{\prime\prime}(h),g_{2}^{\prime\prime}(h))+2(g_{1}^{\prime}(h),g_{2}^{\prime}(h))(A_{2}h+B_{2})^{-1}+(g_{1}(h),g_{2}(h))\\
&\times(A_{2}h+B_{2})^{-1}(I-A_{2})(A_{2}h+B_{2})^{-1})V_{2}(h)\\
=&\frac{1}{G^{2}_{2}(h)}\left[g_{31}(h)I_{11}(h)+g_{32}(h)I_{13}(h)\right].
\end{aligned}
$$
Hence (3.13) holds. By lemma 3.4, we have
$$
\tilde{\Theta}_{2}(h)\triangleq L_{1}(h)\Theta_{2}(h)=\frac{1}{G^{2}_{2}(h)}[Q_{1}(h)I_{11}+Q_{2}(h)I_{13}].
$$
By the same discussion as above, we get
$$
\tilde{\Theta}_{1}(h)\triangleq L_{1}(h)\Theta_{1}(h)=\frac{1}{G_{1}(h)}[Q_{3}(h)I^{\prime}_{01}+Q_{4}(h)Z^{\prime}]
+Q_{5}(h)I^{\prime}_{03}+Q_{6}(h)I^{\prime}_{21}+Q_{7}(h)I^{\prime}_{12}.
$$
Then we have
$$F_{1}(h)=L_{1}(h)I(h)=L_{1}(h)(\Theta_{1}(h)+\Theta_{2}(h))=\tilde{\Theta}_{1}(h)+\tilde{\Theta}_{2}(h).$$
This ends the proof. $\diamondsuit$

Similar to the proof of lemma 3.4 and lemma 3.5, we can get the following lemma,

\vskip 0.2 true cm
{\it \noindent{\bf Lemma 3.6.}~For $h\in\Sigma_{1}$, we obtain

\noindent(i)~There exist polynomials $R_{6}(h)$, $R_{5}(h)$ and $R_{4}(h)$ of $h$ with degree respectively $n_{2}$, $n_{2}-1$ and $n_{2}-2$ such that $L_{2}(h)\tilde{\Theta}_{2}(h)=0$, where $n_{2}=3\left[\frac{n-2}{4}\right]+3\left[\frac{n}{4}\right]+14$, and
$$L_{2}(h)=R_{6}(h)\frac{d^{2}}{dh^{2}}+R_{5}(h)\frac{d}{dh}+R_{4}(h).$$

\noindent(ii)~$L_{2}(h)F_{1}(h)=L_{2}(h)\tilde{\Theta}_{1}(h)\triangleq F_{2}(h)$, where
$$F_{2}(h)=\frac{1}{G^{3}_{1}(h)}[\tilde{Q}_{3}(h)I^{\prime}_{01}+\tilde{Q}_{4}(h)Z^{\prime}]
+\tilde{Q}_{5}(h)I^{\prime}_{03}+\tilde{Q}_{6}(h)I^{\prime}_{21}+\tilde{Q}_{7}(h)I^{\prime}_{12},$$
where
$$\begin{aligned}
{\rm deg} \tilde{Q}_{3}(h)\leq n_{1}&+n_{2}+\left[\frac{n-1}{4}\right]+9, {\rm deg} \tilde{Q}_{4}(h)\leq n_{1}+n_{2}+\left[\frac{n-1}{4}\right]+8,\\
{\rm deg}& \tilde{Q}_{5}(h),\tilde{Q}_{6}(h),\tilde{Q}_{7}(h)\leq n_{1}+n_{2}+\left[\frac{n-3}{4}\right]-3.
\end{aligned}$$
}

\vskip 0.2 true cm
{\it \noindent{\bf Lemma 3.7.}~(i)~$\Theta_{3}(h)$ has at most $\left[\frac{n-2}{4}\right]+2\left[\frac{n}{4}\right]$ zeros on $h\in\Sigma_{1}$, taking into account the multiplicity.

\noindent~(ii)~$\tilde{\Theta}_{2}(h)$ has at most $4\left[\frac{n-2}{4}\right]+5\left[\frac{n}{4}\right]+12$ zeros on $h\in\Sigma_{1}$, taking into account the multiplicity.
}

\noindent{\bf Proof.}~Since $I_{02}\neq0$ on $h\in\Sigma_1$, we have
$$\Theta_{3}(h)=l_{1}(h)I_{02}(h)+l_{2}(h)I_{22}(h)=I_{02}(h)(l_{1}(h)+l_{2}(h)\omega_{3}(h)).$$
Let $W_{3}(h)=l_{1}(h)+l_{2}(h)\omega_{3}(h)$, then $\Theta_{3}(h)$ has the same zeros as $W_{3}(h)$. In the following, we will estimate the number of zeros of $W_{3}(h)$.

By the Riccati equation in lemma 3.3, we have
$$G_{3}(h)l_{2}(h)W_{3}^{\prime}(h)=-b_{2}(h)W_{3}^{2}(h)+m_{1}(h)W_{3}(h)+m_{2}(h),$$
where
\begin{align*}
m_{1}(h)=&2b_{2}(h)l_{1}(h)+(b_{4}(h)-b_{1}(h))l_{2}(h)+G_{3}(h)l_{2}^{\prime}(h),\\ m_{2}(h)=&-l_{1}(h)(2b_{2}(h)l_{1}(h)+(b_{4}(h)-b_{1}(h))l_{2}(h)+G_{3}(h)l_{2}^{\prime}(h))\\
&+b_{2}(h)l^{2}_{1}(h)+l_{2}^{2}(h)b_{3}(h)+G_{3}(h)l_{2}(h)l_{1}^{\prime}(h).
\end{align*}
Then by lemma 4.4 in \cite{zhaoyulin1999}, we get
$$
\begin{aligned}
\#\{\Theta_{3}(h)=0, h\in\Sigma_{1}\}&=\#\{W_{3}(h)=0, h\in\Sigma_{1}\}\\
&\leq\#\{l_{2}(h)=0, h\in\Sigma_{1}\}+\#\{m_{2}(h)=0, h\in\Sigma_{1}\}+1.
\end{aligned}$$
From lemma 2.2, we can get the conclusion (i). The conclusion (ii) can be proved similarly. This ends the proof. $\diamondsuit$

\vskip 0.2 true cm
{\it \noindent{\bf Lemma 3.8.}~For case $D^{+}_{6}$, $I(h)$ has at most $208n+1089$ zeros for $h\in\Sigma_{1}\cup\Sigma_{2}\cup\Sigma_{3}$ (taking into account the multiplicity), where $\Sigma_{i},\ i=1,2,3$ are defined in (3.1).}

\noindent{\bf Proof.}~Recall that system (1.4) has seven singular points:
\\[10pt]
\noindent$\bullet$ Four elementary centers $C_{1,2,3,4}(\pm\sqrt{\frac{b+2c}{b^2-4ac}},\pm\sqrt{-\frac{b+2a}{b^2-4ac}})$,
\vspace{0.5em}

\noindent$\bullet$ Three elementary saddles $S_{0}(0,0)$, $S_{1,2}(0,\pm\sqrt{\frac{1}{2c}})$.

The level curve $H(x,y)=\frac{a+b+c}{b^2-4ac}$ corresponds to the centers $C_{i}$ $(i=1,2,3,4)$, the level curve $H(x,y)=-\frac{1}{4c}$ corresponds to two saddle points $S_{1,2}$, the level curve $H(x,y)=0$ corresponds to saddle $S_{0}$. The system (1.4) has seven period annuli. The period annulus only around $C_{i}$ $(i=1,2,3,4)$ corresponds to $h\in\Sigma_{1}$, the period annulus expanding from the double-homoclinic loop of $S_{1,2}$ to the double-homoclinic loop of $S_{0}$ corresponds $h\in\Sigma_{2}$, the period annulus expanding outwards from the double-homoclinic loop of $S_{0}$ corresponds to $h\in\Sigma_{3}$.

{\bf (1).Zeros of} $I(h)$ {\bf for} $h\in\Sigma_{1}$

Noting that $F_{1}(h)=\tilde{\Theta}_{1}(h)+\tilde{\Theta}_{2}(h)$ and $D(h)\neq0$ for $h\in\Sigma_{1}$, from lemma 3.6, we have $L_{2}(h)F_{1}(h)=F_{2}(h)$. From (3.9), we have
$$
\begin{aligned}
(F_{2}(h))^{(k)}=&\frac{1}{(G_{1}(h))^{k+3}}[L_{k1}(h)I_{01}^{\prime}(h)+L_{k2}(h)Z^{\prime}(h)]\\
=&\frac{I_{01}^{\prime}(h)}{(G_{1}(h))^{k+3}}[L_{k1}(h)+L_{k2}(h)\frac{Z^{\prime}(h)}{I_{01}^{\prime}(h)}],
\end{aligned}
$$
where $k=n_{1}+n_{2}+\left[\frac{n-3}{4}\right]-2$, and
$$
\begin{aligned}
{\rm deg} L_{k1}(h)&\leq n_{1}+n_{2}+\left[\frac{n-1}{4}\right]+9+3k,\\
{\rm deg} L_{k2}(h)&\leq n_{1}+n_{2}+\left[\frac{n-1}{4}\right]+8+3k.
\end{aligned}
$$

Let $W_{1}(h)=L_{k1}(h)+L_{k2}(h)\omega_{1}(h)$, Following the lines of lemma 3.7, and by lemma 3.3, we have
$$\#\{F_{2}(h)=0,h\in\Sigma_{1}\}\leq3\left[\frac{n-1}{4}\right]+10\left[\frac{n-3}{4}\right]+52\left[\frac{n-2}{4}\right]+52\left[\frac{n}{4}\right]+217$$
for $h\in\Sigma_{1}$. Then we have
\begin{align*}
\#\{F_{1}(h)=0, h\in\Sigma_{1}\}\leq&3\#\{R_{6}(h)=0, h\in\Sigma_{1}\}+3\#\{\Theta_{2}(h)=0, h\in\Sigma_{1}\}\\&+\#\{F_{2}(h)=0,h\in\Sigma_{1}\}+2\\
\leq&10\left[\frac{n-3}{4}\right]+73\left[\frac{n-2}{4}\right]+76\left[\frac{n}{4}\right]+3\left[\frac{n-1}{4}\right]+295.
\end{align*}
Noting that there are four period annuli for $h\in\Sigma_{1}$, we obtain from proposition 3.3 that
\begin{align*}
\#\{I(h)=0, h\in\Sigma_{1}\}\leq&4(3\#\{R_{3}(h)=0, h\in\Sigma_{1}\}+3\#\{\Theta_{2}(h)=0, h\in\Sigma_{1}\}\\&+\#\{F_{1}(h)=0\}+2)\\
\leq&40\left[\frac{n-3}{4}\right]+316\left[\frac{n-2}{4}\right]+340\left[\frac{n}{4}\right]+12\left[\frac{n-1}{4}\right]+1212.
\end{align*}

{\bf (2).Zeros of} $I(h)$ {\bf for} $h\in\Sigma_{2}$

By the symmetry, we obtain $I_{11}(h)=I_{12}(h)=I_{13}(h)=0$ for $h\in\Sigma_{2}$. Hence, for $h\in\Sigma_{2}$,
$$I(h)=f_{1}(h)I_{01}+f_{2}(h)I_{03}+f_{3}(h)I_{21}+f_{4}(h)I_{23}+l_{1}(h)I_{02}+l_{2}(h)I_{22}.$$
By Proposition 3.3, we have
$$\begin{aligned}\#\{I(h)=0,h\in\Sigma_{2}\}\leq&2(3\#\{R_{3}(h)=0, h\in\Sigma_{2}\}+3\#\{\Theta_{3}(h)=0, h\in\Sigma_{2}\}\\
&+\#\{F_{3}(h)=0, h\in\Sigma_{2}\}+2),\end{aligned}$$
where
$$F_{3}(h)=L_{1}(h)\Theta_{1}(h)=R_{3}(h)\frac{d^{2}\Theta_{1}(h)}{dh^2}+R_{2}(h)\frac{d\Theta_{1}(h)}{dh}+R_{1}(h)\Theta_{1}(h).$$
Noting that the zero $h=-\frac{1}{4a}$ of $G_{1}(h)$ divides the interval $\Sigma_{2}$ into $\Sigma_{21}=(-\frac{1}{4c},-\frac{1}{4a})$ and $\Sigma_{22}=(-\frac{1}{4a},0)$, then on $\Sigma_{21}\cup\Sigma_{22}$ we have
$$F_{3}(h)=\frac{1}{G_{1}(h)}(P_{0}(h)I^{\prime}_{01}+P_{1}(h)Z^{\prime})+P_{2}(h)I_{03}^{\prime}+P_{3}(h)I_{21}^{\prime}$$
where
$$\begin{aligned}
{\rm deg}P_{0}(h)&\leq\left[\frac{n-1}{4}\right]+\left[\frac{n-2}{4}\right]+\left[\frac{n}{4}\right]+5,\quad {\rm deg}P_{1}(h)\leq\left[\frac{n-1}{4}\right]+\left[\frac{n-2}{4}\right]+\left[\frac{n}{4}\right]+4,\\
{\rm deg}P_{2}(h)&\leq\left[\frac{n-3}{4}\right]+\left[\frac{n-2}{4}\right]+\left[\frac{n}{4}\right]+1,\quad {\rm deg}P_{3}(h)\leq\left[\frac{n-3}{4}\right]+\left[\frac{n-2}{4}\right]+\left[\frac{n}{4}\right]+1.
\end{aligned}$$

Similar to the discussion of the number of zeros of $F_{2}(h)$, we get
$$\#\{F_{3}(h)=0, h\in\Sigma_{21}\cup\Sigma_{22}\}\leq 10\left[\frac{n-3}{4}\right]+13\left[\frac{n-2}{4}\right]+13\left[\frac{n}{4}\right]+3\left[\frac{n-1}{4}\right]+37,\\
$$
and
$$\#\{F_{3}(h)=0, h\in\Sigma_{2}\}\leq\#\{F_{3}(h)=0, h\in\Sigma_{21}\cup\Sigma_{22}\}+1.$$
Hence we obtain
$$\#\{I(h)=0,h\in\Sigma_{2}\}\leq20\left[\frac{n-3}{4}\right]+38\left[\frac{n-2}{4}\right]+44\left[\frac{n}{4}\right]+6\left[\frac{n-1}{4}\right]+92.$$

{\bf (3).Zeros of} $I(h)$ {\bf for} $h\in\Sigma_{3}$

By the symmetry, we obtain $I_{12}(h)=I_{11}(h)=I_{13}(h)=I_{02}(h)=I_{22}(h)=0$ for $h\in\Sigma_{3}$. Hence, for $h\in\Sigma_{3}$,
$$I(h)=f_{1}(h)I_{01}+f_{2}(h)I_{03}+f_{3}(h)I_{21}+f_{4}(h)I_{23}.$$
We split the proof into three steps.

(i)
From (3.9), we obtain
$$
\begin{aligned}
I(h)&=u_{1}(h)I_{01}^{\prime}+u_{2}(h)I_{03}^{\prime}+u_{3}(h)I_{21}^{\prime}+u_{4}(h)Z^{\prime},\\
I^{\prime}(h)&=v_{1}(h)I_{01}^{\prime}+v_{2}(h)I_{03}^{\prime}+v_{3}(h)I_{21}^{\prime}+v_{4}(h)Z^{\prime},
\end{aligned}
\eqno(3.14)$$
where $u_{i}(h)$ and $v_{i}(h)$ for $i=1,2,3,4$ are polynomials of $h$ with
$$
\begin{aligned}
{\rm deg}& u_{1}(h)\leq\left[\frac{n-1}{4}\right]+1, {\rm deg} u_{2}(h),u_{3}(h)\leq\left[\frac{n-3}{4}\right]+1, {\rm deg} u_{4}(h)\leq\left[\frac{n-1}{4}\right],\\
&{\rm deg} v_{1}(h)\leq\left[\frac{n-1}{4}\right], {\rm deg} v_{2}(h),v_{3}(h)\leq\left[\frac{n-3}{4}\right], {\rm deg} v_{4}(h)\leq\left[\frac{n-1}{4}\right]-1.
\end{aligned}
$$
Removing $I_{21}^{\prime}$ from (3.14), we get
$$u_{3}(h)I^{\prime}(h)=v_{3}(h)I(h)+C(h),$$
where
$$C(h)=c_{1}(h)I_{01}^{\prime}+c_{2}(h)I_{03}^{\prime}+c_{3}(h)Z^{\prime}$$
with deg$c_{1}(h)\leq\left[\frac{n-1}{4}\right]+\left[\frac{n-3}{4}\right]+1$, deg$c_{2}(h)\leq2\left[\frac{n-3}{4}\right]+1$ and deg$c_{3}(h)\leq\left[\frac{n-1}{4}\right]+\left[\frac{n-3}{4}\right].$
It follows form lemma 5.1 in \cite{zhangzhifen2002} that
$$\#\{I(h)=0, h\in\Sigma_{3}\}\leq\#\{u_{3}(h)=0,h\in\Sigma_{2}\}+\#\{C(h)=0,h\in\Sigma_{2}\}+1.\eqno(3.15)$$

(ii)
For $h\in\Sigma_{3}$, $G_{1}(h)\neq0$. Derivating $C(h)$ with respect to $h$ gives
$$G_{1}(h)c_{2}(h)C^{\prime}(h)=G_{1}(h)c^{\prime}_{2}(h)C(h)+C_{1}(h),$$
where
$$
\begin{aligned}
C_{1}(h)=&\left[
c_{1}^{\prime}(h)c_{2}(h)G_{1}(h)+d_{11}(h)c_{1}(h)c_{2}(h)+d_{21}c_{2}^{2}(h)+d_{41}(h)c_{2}(h)c_{3}(h)\right.\notag\\
&\left.-c_{1}(h)c_{2}^{\prime}(h)G_{1}(h)\right]I_{01}^{\prime}+\left[c_{2}(h)c_{3}^{\prime}(h)G_{1}(h)+d_{12}(h)c_{1}(h)c_{2}(h)\right.\notag\\
&\left.+d_{22}(h)c_{2}^{2}(h)+d_{42}(h)c_{2}(h)c_{3}(h)-c_{2}^{\prime}(h)c_{3}(h)G_{1}(h)\right]Z^{\prime}\notag\\
=&c_{11}(h)I_{01}^{\prime}+c_{12}(h)Z^{\prime}
\end{aligned}\eqno(3.16)
$$
with deg$c_{11}(h)\leq3\left[\frac{n-3}{4}\right]+\left[\frac{n-1}{4}\right]+5$, deg$c_{12}(h)\leq3\left[\frac{n-3}{4}\right]+\left[\frac{n-1}{4}\right]+4.$
By lemma 5.1 in \cite{zhangzhifen2002}, we get
$$\#\{C(h)=0,h\in\Sigma_{3}\}\leq\#\{c_{2}(h)=0,h\in\Sigma_{3}\}+\#\{C_{1}(h)=0,h\in\Sigma_{3}\}+1.\eqno(3.17)$$

(iii)
From (3.16), we get
$$C_{1}(h)=I_{01}^{\prime}\left[c_{11}(h)+c_{12}(h)\omega_{1}(h)\right].$$
Let $W_{3}(h)=c_{11}(h)+c_{12}(h)\omega_{1}(h)$, then by lemma 3.3 we get
$$G_{1}(h)c_{12}(h)W_{3}^{\prime}(h)=-d_{12}(h)W_{3}^{2}(h)+C_{2}(h)W_{3}(h)+C_{3}(h),$$
where \begin{align*}
C_{3}(h)=&-(G_{1}(h)c_{11}(h)c_{12}^{\prime}(h)+G_{1}(h)c_{11}^{\prime}(h)c_{12}(h)-d_{11}(h)c_{11}(h)c_{12}(h)+d_{12}(h)c_{11}^{2}(h)\\
&-d_{41}(h)c_{12}^{2}(h)+d_{42}(h)c_{11}(h)c_{12}(h))
\end{align*}
with deg$C_{3}(h)\leq2\left[\frac{n-1}{4}\right]+6\left[\frac{n-3}{4}\right]+12$. By lemma 4.4 in \cite{zhaoyulin1999}, for $h\in\Sigma_{3}$, we obtain
$$
\#\{C_{1}(h)=0,h\in\Sigma_{3}\}\leq\#\{c_{12}(h)=0,h\in\Sigma_{3}\}+\#\{C_{3}(h)=0,h\in\Sigma_{3}\}+1.
\eqno(3.18)$$
From (3.15), (3.17) and (3.18), we obtain
$$
\#\{I(h)=0, h\in\Sigma_{3}\}\leq12\left[\frac{n-3}{4}\right]+3\left[\frac{n-1}{4}\right]+21.
$$
Combined {\bf(1)} and {\bf(2)}, we have
$$
\begin{aligned}
\#\{I(h)=0, h\in\Sigma_{1}\cup\Sigma_{2}\cup\Sigma_{3}\}\leq&21\left[\frac{n-1}{4}\right]+354\left[\frac{n-2}{4}\right]+72\left[\frac{n-3}{4}\right]
+384\left[\frac{n}{4}\right]+1325\\
<&208n+1089.
\end{aligned}
$$
This ends the proof. $\diamondsuit$

Let $\Sigma$ is an union of several open intervals on which the period annuli are defined for each case in Fig 1 and Fig 2, we obtain the following three lemmas.

\vskip 0.2 true cm
{\it \noindent{\bf Lemma 3.9.}~For $ac(b^2-4ac)\neq0$, $\#\{G_{1}(h)=0, h\in\Sigma\}=0$, we have

\noindent(i)~$\#\{I(h)=0, (a,b,c)\in D_{1(1)}^{+}\cup D_{1(3)}^{+}\}\leq58n+121$,

\noindent(ii)~$\#\{I(h)=0, (a,b,c)\in D_{1(2)}^{+}\}\leq54n+109$,

\noindent(iii)~$\#\{I(h)=0, (a,b,c)\in l_{1}^{+}\cup l_{2}^{+}\cup D_{4(1)}^{+}\cup D_{1(3)}^{-}\}\leq31n+66$,

\noindent(iv)~$\#\{I(h)=0, (a,b,c)\in D_{4(2)}^{+}\cup D_{4(3)}^{+}\cup D_{1(1)}^{-}\cup D_{1(2)}^{-}\}\leq27n+55$.
}

\noindent{\bf Proof.}~By the same arguments as the proof of lemma 3.8, we get conclusions (i)--(iv). This ends the proof. $\diamondsuit$

\vskip 0.2 true cm
{\it \noindent{\bf Lemma 3.10}~For $ac(b^2-4ac)\neq0$, $\#\{G_{1}(h)=0, h\in\Sigma\}\neq0$, we have

\noindent(i)~$\#\{I(h)=0, (a,b,c)\in D_{2}^{+}\cup D_{3}^{+}\}\leq31n+67$,

\noindent(ii)~$\#\{I(h)=0, (a,b,c)\in l_{3}^{+}\cup l_{1}^{-}\}\leq31n+68$,

\noindent(iii)~$\#\{I(h)=0, (a,b,c)\in D_{5}^{+}\cup D_{2}^{-}\}\leq31n+69$,

\noindent(iv)~$\#\{I(h)=0, (a,b,c)\in D_{6}^{+}\cup D_{3}^{-}\}\leq208n+1089$.
}

\noindent{\bf Proof.}~Noting that $G_{1}(h)$ has four zeros
$$h_{1}=0,\ h_{2}=-\frac{1}{4c},\ h_{3}=-\frac{1}{4a},\ h_{4}=\frac{a+b+c}{b^2-4ac},$$
the position of $h_{i},i=1,2,3,4$ are summarized in Table 2 for each case in $S$. By the same arguments as the proof of lemma 3.8, we get conclusions (i)--(iv). This ends the proof. $\diamondsuit$
\begin{table}[htbp]
\centering
\resizebox{\linewidth}{!}{\scriptsize
\begin{tabular}{cccc}
 \hline
  \makebox[0.2\textwidth][c]{\multirow{2}{*}{Case}}&\makebox[0.2\textwidth][c]{\multirow{2}{*}{Parameter range}}&\makebox[0.2\textwidth][c]{\multirow{2}{*}{Period annuli}}&\makebox[0.2\textwidth][c]{\multirow{2}{*}{Position of zeros}} \\[10pt]
   \hline
\multirow{2}{*}{$D^{+}_{2}$}    &\multirow{2}{*}{-}&\multirow{2}{*}{$(-\frac{1}{4c},0)\cup(0,-\frac{1}{4a})$} &\multirow{2}{*}{$-\frac{1}{4c}<0<h_{4}<\frac{1}{4a}$}  \\[10pt]
 \multirow{2}{*}{$D^{+}_{3}$} &\multirow{2}{*}{-}&\multirow{2}{*}{$(-\frac{1}{4c},0)\cup(0,-\frac{1}{4a})$}  &\multirow{2}{*}{$-\frac{1}{4c}<h_{4}<0<\frac{1}{4a}$} \\[10pt]
 \multirow{2}{*}{$D^{+}_{5}$}    &\multirow{2}{*}{$b^2-4ac<0,a<c$} &\multirow{2}{*}{$(-\frac{1}{4c},0)\cup(0,+\infty)$} &\multirow{2}{*}{$h_{4}<h_{3}<-\frac{1}{4c}<0$} \\[10pt]
 \multirow{2}{*}{$D^{+}_{5}$} &\multirow{2}{*}{$b^2-4ac<0,a>c$}&\multirow{2}{*}{$(-\frac{1}{4c},0)\cup(0,+\infty)$}  &\multirow{2}{*}{$h_{4}<-\frac{1}{4c}<h_{3}<0$} \\[10pt]
  \multirow{2}{*}{$D^{+}_{5}$}    &\multirow{2}{*}{$b^2-4ac>0,a<c$}&\multirow{2}{*}{$(-\frac{1}{4c},0)\cup(0,+\infty)$} &\multirow{2}{*}{$h_{3}<-\frac{1}{4c}<0<h_{4}$} \\[10pt]
   \multirow{2}{*}{$D^{+}_{5}$}    &\multirow{2}{*}{$b^2-4ac>0,a>c$}&\multirow{2}{*}{$(-\frac{1}{4c},0)\cup(0,+\infty)$} &\multirow{2}{*}{$-\frac{1}{4c}<h_{3}<0<h_{4}$} \\[10pt]
    \multirow{2}{*}{$l^{+}_{3}$}    &\multirow{2}{*}{-} &\multirow{2}{*}{$(-\frac{1}{4c},0)\cup(0,+\infty)$}  &\multirow{2}{*}{$h_{4}<-\frac{1}{4c}<h_{3}<0$} \\[10pt]
    \multirow{2}{*}{$D^{+}_{6}$}    &\multirow{2}{*}{-} &\multirow{2}{*}{$(\frac{a+b+c}{b^2-4ac},-\frac{1}{4c})\cup(-\frac{1}{4c},0)\cup(0,+\infty)$}  &\multirow{2}{*}{$\frac{a+b+c}{b^2-4ac}<-\frac{1}{4c}<h_{3}<0$} \\[10pt]
    \multirow{2}{*}{$D^{-}_{2}$}    &\multirow{2}{*}{$b^2-4ac<0,a<c$} &\multirow{2}{*}{$(-\infty,0)\cup(0,-\frac{1}{4a})$}  &\multirow{2}{*}{$0<-\frac{1}{4a}<h_{2}<h_{4}$} \\[10pt]
    \multirow{2}{*}{$D^{-}_{2}$}    &\multirow{2}{*}{$b^2-4ac<0,a>c$} &\multirow{2}{*}{$(-\infty,0)\cup(0,-\frac{1}{4a})$}  &\multirow{2}{*}{$0<h_{2}<-\frac{1}{4a}<h_{4}$} \\[10pt]
    \multirow{2}{*}{$D^{-}_{2}$}    &\multirow{2}{*}{$b^2-4ac>0,a<c$} &\multirow{2}{*}{$(-\infty,0)\cup(0,-\frac{1}{4a})$}  &\multirow{2}{*}{$h_{4}<0<-\frac{1}{4a}<h_{2}$} \\[10pt]
    \multirow{2}{*}{$D^{-}_{2}$}    &\multirow{2}{*}{$b^2-4ac>0,a>c$} &\multirow{2}{*}{$(-\infty,0)\cup(0,-\frac{1}{4a})$}  &\multirow{2}{*}{$h_{4}<0<h_{2}<-\frac{1}{4a}$} \\[10pt]
    \multirow{2}{*}{$l^{-}_{1}$}    &\multirow{2}{*}{-} &\multirow{2}{*}{$(-\infty,0)\cup(0,-\frac{1}{4a})$}  &\multirow{2}{*}{$0<h_{2}<-\frac{1}{4a}<h_{4}$} \\[10pt]
    \multirow{2}{*}{$D^{-}_{3}$}    &\multirow{2}{*}{-} &\multirow{2}{*}{$(-\infty,0)\cup(0,-\frac{1}{4a})\cup(-\frac{1}{4a},\frac{a+b+c}{b^2-4ac})$}  &\multirow{2}{*}{$0<h_{2}<-\frac{1}{4a}<\frac{a+b+c}{b^2-4ac}$} \\[10pt]
   \hline
 \end{tabular}}
\multirow{4}{*}{Table 2. The position of zeros of $G_{1}(h)$ for each case in $S$}
 \end{table}

\vskip 0.2 true cm
{\it \noindent{\bf Lemma 3.11}~For $a=0, c>0$, we have

\noindent(i)~$\#\{I(h)=0, (a,b,c)\in D_{4(1)}^{+}\cup D_{5}^{+}\}\leq49n-60$,

\noindent(ii)~$\#\{I(h)=0, (a,b,c)\in D_{4(2)}^{+}\cup D_{4(3)}^{+}\}\leq45n-58$,

\noindent(iii)~$\#\{I(h)=0, a=b=0,c>0\}\leq9n-1$.
}

\noindent{\bf Proof.}~(1)~For $D_{4(1)}^{+},D_{4(2)}^{+},D_{4(3)}^{+}$, and $D_{5}^{+}$ with $a=0$, $b\neq0$, we consider an equivalent system of (1.5) whose unperturbed system has the Hamiltonian
$$H(x,y)=y^2-x^2+bx^2y^2+cx^4.$$
The corresponding $I(h)$ can be expressed as
$$I(h)=\alpha(h)I_{01}(h)+\beta(h)I_{03}(h)+\gamma(h)I_{21}(h)+\phi(h)I_{11}(h)+\psi(h)I_{13}(h),$$
and the following Picard-Fuchs equations satisfied by $V_{5}=(I_{01},I_{03},I_{21})^{T}$ and $V_{6}=(I_{11},I_{13})^{T}$,
$$V_{5}=(A_{5}h+B_{5})V_{5}^{\prime},\quad V_{6}=(A_{6}h+B_{6})V_{6}^{\prime}.$$
where $\alpha(h),\beta(h),\gamma(h),\phi(h),\psi(h)$ are real polynomials in $h$ with
$${\rm deg}\alpha(h),\beta(h),\gamma(h)\leq\left[\frac{n-1}{2}\right]-1,\quad{\rm deg}\phi(h),\psi(h)\leq\left[\frac{n}{2}\right]-2,$$
and
$$A_5=\left( \begin {array}{cccc} 2&0&0\\ \noalign{\medskip}\frac{5c+4b}{2b^2}
&\frac{2}{3}&-\frac{c}{b}\\ \noalign{\medskip}-\frac{1}{2b}&0&1\end {array}
 \right),\qquad
B_5=\left( \begin {array}{cccc} 0&-\frac{1}{3}&1\\ \noalign{\medskip}0&-\frac{13b+15c}{12b^2}&\frac{3b+5c}{4b^2}\\ \noalign{\medskip}0&\frac{3c+b}{12bc}&\frac{b-c}{4bc}\end {array}
 \right),$$
$$A_6=\left( \begin {array}{cccc} \frac43&0\\ \noalign{\medskip}\frac{4bc+8c^2}{5b^2c}&\frac45\end {array}
 \right),\qquad
B_6=\left( \begin {array}{cccc} \frac{1}{3c}&-\frac{b+2c}{9c}\\  \noalign{\medskip}\frac{b+2c}{5b^2c}&-\frac{b^2+16bc+16c^2}{15b^2c}\end {array}
 \right),$$
Taking
$$\bar{Z}=\frac{2}{5}I_{03}+\frac{6c}{5b}I_{21},$$
we obtain
$$
G_{4}(h) \left(
\begin{matrix}
I_{01}^{\prime\prime}\\
I_{03}^{\prime\prime}\\
\bar{Z}^{\prime\prime}\\
\end{matrix}
\right)=\left(
\begin{matrix}
\bar{d}_{11}(h) \quad \bar{d}_{12}(h)\\
\bar{d}_{21}(h) \quad \bar{d}_{22}(h)\\
\bar{d}_{31}(h) \quad \bar{d}_{32}(h)\\
\end{matrix}
\right)\left(
\begin{matrix}
I_{01}^{\prime}\\
\bar{Z}^{\prime}\\
\end{matrix}
\right),
$$
where $G_{4}(h)=b^2h(4ch+1)(h-\frac{b+c}{b^2})$, and
\begin{equation*}
\begin{tabular}{lll}
$\bar{d}_{11}(h)=-\frac12h(4b^2ch+b^2-bc-2c^2),$&\qquad\qquad&$\bar{d}_{12}(h)=\frac{5b}{12}(2bch+b+c),$\\
$\ $&\qquad\qquad&$\ $\\
$\bar{d}_{21}(h)=-\frac32(b+c)h(4ch+1),$&\qquad\qquad&$\bar{d}_{22}(h)=\frac{5b^2}{4}h(4ch+1),$\\
$\ $&\qquad\qquad&$\ $\\
$\bar{d}_{31}(h)=-\frac{3}{5b}h(4b^2ch+2bc^2h+b^2-c^2),$&\qquad\qquad&$\bar{d}_{32}(h)=\frac12h(4b^2ch+b^2-bc-2c^2).$
\end{tabular}
\end{equation*}

For $I_{11}$ and $I_{13}$, we have
$$
G_{5}(h) \left(
\begin{matrix}
I_{11}^{\prime}\\
I_{13}^{\prime}\\
\end{matrix}
\right)=\left(
\begin{matrix}
\bar{a}_{1}(h) \quad \bar{a}_{2}(h)\\
\bar{a}_{3}(h) \quad \bar{a}_{4}(h)\\
\end{matrix}
\right)\left(
\begin{matrix}
I_{11}\\
I_{13}\\
\end{matrix}
\right),
$$
where $G_{5}(h)=b^2(4ch+1)(h-\frac{b+c}{b^2})$ and

\begin{equation*}
\begin{tabular}{lll}
$\bar{a}_{1}(h)=\frac{1}{4}(12b^2ch-b^2-16bc-16c^2),$&\qquad\qquad&$\bar{a}_{2}(h)=\frac{5}{12}(b+2c)b^2,$\\
$\ $&\qquad\qquad&$\ $\\
$\bar{a}_{3}(h)=-\frac{3}{4}(b+2c)(4ch+1),$&\qquad\qquad&$\bar{a}_{4}(h)=\frac{5}{4}b^2(4ch+1).$
\end{tabular}
\end{equation*}
Set $\bar{\omega}_{1}(h)=\frac{\bar{Z}}{I_{01}}$, $\bar{\omega}_{2}(h)=\frac{I_{13}}{I_{11}}$, similar to lemma 3.3, we get the following two Riccati equations
$$G_{4}(h){\bar\omega}_{1}^{\prime}(h)=-\bar{d}_{12}{\bar\omega}_{1}^{2}(h)+(\bar{d}_{32}(h)-\bar{d}_{11}(h)){\bar\omega}_{1}(h)+\bar{d}_{31}(h),$$
$$G_{5}(h){\bar\omega}_{2}^{\prime}(h)=-\bar{a}_{2}{\bar\omega}_{2}^{2}(h)+(\bar{a}_{4}(h)-\bar{a}_{1}(h)){\bar\omega}_{2}(h)+\bar{a}_{3}(h).$$
Noting that $I(h)=CI_{01}(h)$ when $n=1$, where $C$ is a arbitrary constant. If $C\neq0$, $I(h)$ has no zeros as $I_{01}(h)\neq0$ on $(-\frac{1}{4c},0)\cup(0,+\infty)$. Following the lines of lemma 3.8, we can get the conclusion (i) and (ii) for $n\geq2$.

\noindent(2)~For $D_{5}^{+}$ with $a=b=0$, $c>0$, we consider an equivalent system of (1.5) whose unperturbed system has the Hamiltonian
$$H(x,y)=y^2-x^2+cx^4.$$
The corresponding $I(h)$ can be expressed as
$$I(h)=\sigma(h)I_{01}(h)+\delta(h)I_{11}(h)+\zeta(h)I_{21}(h),$$
where $\sigma(h),\delta(h),\zeta(h)$ are real polynomials in $h$ with
$${\rm deg}\sigma(h)\leq\left[\frac{n-1}{2}\right],\quad{\rm deg}\delta(h)\leq\left[\frac{n}{2}\right]-1,\quad{\rm deg}\zeta(h)\leq\left[\frac{n-1}{2}\right]-1.$$
For $I_{01}$ and $I_{21}$, we obtain
$$
h(4ch+1) \left(
\begin{matrix}
I_{01}^{\prime}\\
I_{21}^{\prime}\\
\end{matrix}
\right)=\left(
\begin{array}{cc}
3ch+1&-\frac{5c}{2}\\
-\frac{h}{2}&5ch\\
\end{array}
\right)\left(
\begin{matrix}
I_{01}\\
I_{21}\\
\end{matrix}
\right).
$$
Generator $I_{11}$ satisfy the following homogeneous linear differential equation
$$I_{11}=(h+\frac{1}{4c})I_{11}^{\prime},$$
we get $I_{11}=C_{1}(h+\frac{1}{4c})$, where $C_{1}=4cI_{11}(0)$. Then we obtain
$$I(h)=\sigma(h)I_{01}(h)+\zeta(h)I_{21}(h)+\tilde{\delta}(h),$$
where ${\rm deg}\tilde{\delta}(h)\leq\left[\frac{n}{2}\right]$.
Following the lines of lemma 3.8, we can get the conclusion (iii). This ends the proof. $\diamondsuit$

{\bf\noindent Proof of Theorem 1.1}~Following the line of lemma 3.9, lemma 3.10 and lemma 3.11, we can obtain Theorem 1.1 for $n>4$. We can also prove that Theorem 1.1 is true for $n\leq4$. This ends the proof. $\diamondsuit$

\section{Hopf and homoclinic bifurcation of $l_{2}^{+}$}

In this section, we study the Hopf bifurcation and homoclinic bifurcation of (1.6). Then we will give the distributions of limit cycles near the center and double-homoclinic loop.

\vspace{1em}
\noindent{\bf (1)~Hopf bifurcation of $l_{2}^{+}$}

By the symmetry, system (1.6) is equivalent to
$$
\begin{cases}
\dot{x}=2y(1+2x^2-2y^2),\\
\dot{y}=-2x(1-2x^2-2y^2)+\varepsilon(q_{0}y+q_{1}xy+q_{2}y^3+q_{3}x^2y),
\end{cases}
$$
where $q_{0}$, $q_{1}$, $q_{2}$ and $q_{3}$ are arbitrary constants. Then we get
$$
\begin{cases}
\dot{x}=-4y-4\sqrt{2}xy-4x^2y+4y^3,\\
\dot{y}=4x+6\sqrt{2}x^2+4xy^2+4x^3+2\sqrt{2}y^2+\varepsilon(\alpha_{0}y+\alpha_{1}xy+\alpha_{2}y^3+\alpha_{3}x^2y),
\end{cases}\eqno(4.1)
$$
by making variable transformation $x=x+\frac{1}{\sqrt{2}}$, $y=y$, where
$$
\begin{tabular}{lll}
$\alpha_{0}=q_{0}+\frac{\sqrt{2}}{2}q_{1}+\frac{1}{2}q_{3},$&\qquad\qquad&$\alpha_{1}=q_{1}-\sqrt{2}q_{3},$\\
\ &\qquad\qquad&\ \\
$\alpha_{2}=q_{2},$&\qquad\qquad&$\alpha_{3}=q_{3}.$\\
\end{tabular}
$$
with det$\frac{\partial(\alpha_{0},\alpha_{1},\alpha_{2},\alpha_{3})}{\partial(q_{0},q_{1},q_{2},q_{3})}=1$.
System (4.1)$|_{\varepsilon=0}$ has two centers $(0,0)$ and $(-\sqrt{2},0)$. The level curves around two centers correspond to $H_{1}(x,y)=h$, $h\in(-\frac{1}{4},0)$, where
$$H_{1}(x,y)=-2x^2-2y^2-2\sqrt{2}x^3-2\sqrt{2}xy^2-2x^2y^2+y^4-x^4.$$
Then we get the following lemma.

\vskip 0.2 true cm
{\it \noindent{\bf Lemma 4.1.}~For $0<-h\ll1$, $I(h)$ has the following forms near a center $(0,0)$ with $H_{1}(x,y)=0$:
$$I(h)=\sum\limits_{j\geq0}a_{j}(\delta)m^{2j},$$
where $m=\sqrt{-\frac{1}{2}h}$, $a_{j}(\delta)$ are coefficients depending on $\delta=(\alpha_{0},\alpha_{1},\alpha_{2},\alpha_{3})$.
}

\noindent{\bf Proof} By introducing the polar coordinate change $x=rcos\theta$, $y=rsin\theta$ so that the curve $H_{1}(x,y)=h$ can be expressed as
$$r^2(1-\frac{1}{2}r^2+\sqrt{2}rcos\theta+2r^2cos^2\theta-r^2cos^4\theta)=-\frac{1}{2}h.\eqno(4.2)$$
Equation (4.2) has a unique positive solution in $r$, and it can be written as
$$r=m+e_{2}(\theta)m^2+e_{3}(\theta)m^3+\cdots+e_{6}(\theta)m^6+O(m^7)\triangleq\phi(\theta,m),\eqno(4.3)$$
where $m=\sqrt{-\frac{1}{2}h}$, and
$$\begin{aligned}
e_{2}(\theta)=&-\frac{\sqrt{2}}{2}cos\theta,\\
e_{3}(\theta)=&\frac12cos^4\theta+\frac14cos^2\theta+\frac14,\\
e_{4}(\theta)=&-\frac{\sqrt{2}}{4}(6cos^4\theta-4cos^2\theta+3)cos\theta,\\
e_{5}(\theta)=&\frac78cos^8\theta+\frac{35}{8}cos^6\theta-\frac{133}{32}cos^4\theta+\frac{35}{16}cos^2\theta+\frac{7}{32},\\
e_{6}(\theta)=&-\frac{\sqrt{2}}{4}(20cos^8\theta-4cos^4\theta+5)cos\theta.\\
\end{aligned}$$
By Green formula and (4.3), we obtain
$$\begin{aligned}
I(h)=&\alpha_{0}\oint_{\Gamma_{h}}ydx+\alpha_{1}\oint_{\Gamma_{h}}xydx+\alpha_{2}\oint_{\Gamma_{h}}y^{3}dx+\alpha_{3}\oint_{\Gamma_{h}}x^{2}ydx\\
=&\alpha_0\iint_{D}-1dxdy+\alpha_1\iint_{D}-xdxdy+\alpha_2\iint_{D}-3y^2dxdy+\alpha_3\iint_{D}-x^2dxdy\\
=&-\frac{\alpha_0}{2}\int_0^{2\pi}\phi^2(\theta,m)d\theta-\frac{\alpha_1}{3}\int_0^{2\pi}\cos\theta\phi^3(\theta,m)d\theta\\
&-\frac{3\alpha_2}{4}\int_0^{2\pi}\sin^2\theta\phi^4(\theta,m)d\theta-\frac{\alpha_3}{4}\int_0^{2\pi}\cos^2\theta\phi^4(\theta,m)d\theta\notag\\
=&\sum\limits_{j\geq0}a_{j}(\delta)m^{2j},
\end{aligned}$$
where
$$\begin{aligned}
a_{1}(\delta)=&-\alpha_{0}\pi,\\
a_{2}(\delta)=&-(\frac{11}{8}\alpha_{0}\pi-\frac{\sqrt{2}}{2}\alpha_{1}\pi+\frac{3}{4}\alpha_{2}\pi+\frac14\alpha_{3}\pi),\\
a_{3}(\delta)=&-(\frac{259}{64}\alpha_{0}\pi-\frac{7\sqrt{2}}{4}\alpha_{1}\pi+\frac{27}{16}\alpha_{2}\pi+\frac{21}{16}\alpha_{3}\pi),\\
a_{4}(\delta)=&-(\frac{50185}{8192}\alpha_{0}\pi-\frac{1885\sqrt{2}}{256}\alpha_{1}\pi+\frac{2505}{512}\alpha_{2}\pi+\frac{3195}{512}\alpha_{3}\pi),\\
\ldots\\
\end{aligned}$$
This ends the proof. $\diamondsuit$

For simplicity, we consider the number of zeros of
$$M_{1}(h)=-\sum\limits_{j\geq1}a_{j}(\delta)m^{2j}=\sum\limits_{j\geq1}b_{j}(\delta)m^{2j},$$
where $a_{j}(\delta)=-b_{j}(\delta)$. Then we have the following lemma.

\vskip 0.2 true cm
{\it \noindent{\bf Lemma 4.2.}~$M_{1}(h)$ has exactly 3 zeros near $h=0$.

}
\noindent{\bf Proof}
By elementary calculation, we have
$$
\begin{aligned}
&b_{2}|_{b_{1}=0}=-\frac{\sqrt{2}}{2}\alpha_{1}\pi+\frac{3}{4}\alpha_{2}\pi+\frac{1}{4}\alpha_{3}\pi,\\
&b_{3}|_{b_{1}=b_{2}=0}=-\frac{15}{16}\alpha_{2}\pi+\frac{7}{16}\alpha_{3}\pi,\\
&b_{4}|_{b_{1}=b_{2}=b_{3}=0}=-\frac{5}{16}\alpha_{3}\pi,
\end{aligned}
$$
and there exist $\delta_{0}=(0,\frac{3\sqrt{2}}{5}\alpha_{3}^{*},\frac{7}{15}\alpha_{3}^{*},\alpha_{3}^{*})$, where $\alpha_{3}^{*}\neq0$, such that $b_{1}(\delta_{0})=b_{2}(\delta_{0})=b_{3}(\delta_{0})=0$, $b_{4}(\delta_{0})\neq0$. Furthermore, we have
$$\text{det}\frac{\partial(b_{1},b_{2},b_{3})}{\partial(\alpha_{0},\alpha_{1},\alpha_{2})}(\delta)=\frac{15\sqrt{2}}{32}\pi^3\neq0,$$
which yields that $M_{1}(h)$ has exactly 3 zeros near $h=0$. This ends the proof. $\diamondsuit$

\vspace{1em}
For center $(-\sqrt{2},0)$, we get the equivalent system
$$
\begin{cases}
\dot{x}=-4y+4\sqrt{2}xy-4x^2y+4y^3,\\
\dot{y}=4x-6\sqrt{2}x^2+4xy^2+4x^3-2\sqrt{2}y^2+\varepsilon(\hat{\alpha}_{0}y+\hat{\alpha}_{1}xy+\hat{\alpha}_{2}y^3+\hat{\alpha}x^2y),
\end{cases}\eqno(4.4)
$$
of (4.1) by making variable transformation $x=x-\sqrt{2}$, $y=y$. where
$$
\begin{tabular}{lll}
$\hat{\alpha}_{0}=\alpha_{0}-2\sqrt{2}\alpha_{1}+2\alpha_{3},$&\qquad\qquad&$\hat{\alpha}_{1}=\alpha_{1}-2\sqrt{2}\alpha_{3},$\\
\ &\qquad\qquad&\ \\
$\hat{\alpha}_{2}=\alpha_{2},$&\qquad\qquad&$\hat{\alpha}_{3}=\alpha_{3}.$\\
\end{tabular}
$$
with det$\frac{\partial(\hat{\alpha}_{0},\hat{\alpha}_{1},\hat{\alpha}_{2},\hat{\alpha}_{3})}{\partial(\alpha_{0},\alpha_{1},\alpha_{2},\alpha_{3})}=1$. System (4.4)$|_{\varepsilon=0}$ has a center $(0,0)$ corresponding to $H_{2}(0,0)=0$, where
$$H_{2}(x,y)=-2x^2-2y^2+2\sqrt{2}x^3+2\sqrt{2}xy^2-2x^2y^2+y^4-x^4.$$

Similar to the proof of lemma 4.1, we get the following lemma.

\vskip 0.2 true cm
{\it \noindent{\bf Lemma 4.3.}~For $0<-h\ll1$, $I(h)$ has the following forms near a center $(0,0)$ with $H_{2}(0,0)=0$
$$I(h)=-M_{2}(h)=-\sum\limits_{l\geq1}d_{l}(\delta)m^{2l},$$
where
$$
\begin{aligned}
d_{1}(\delta)=&(\alpha_{0}-2\sqrt{2}\alpha_{1}+2\alpha_{3})\pi,\\
d_{2}(\delta)=&-(\frac{11}{8}(\alpha_{0}-2\sqrt{2}\alpha_{1}+2\alpha_{3})\pi+\frac{\sqrt{2}}{2}(\alpha_{1}-2\sqrt{2}\alpha_{3})\pi+\frac{3}{4}\alpha_{2}\pi+\frac14\alpha_{3}\pi),\\
d_{3}(\delta)=&-(\frac{259}{64}(\alpha_{0}-2\sqrt{2}\alpha_{1}+2\alpha_{3})\pi+\frac{7\sqrt{2}}{4}(\alpha_{1}-2\sqrt{2}\alpha_{3})\pi+\frac{27}{16}\alpha_{2}\pi+\frac{21}{16}\alpha_{3}\pi),\\
d_{4}(\delta)=&-(\frac{50185}{8192}(\alpha_{0}-2\sqrt{2}\alpha_{1}+2\alpha_{3})\pi+\frac{1885\sqrt{2}}{256}(\alpha_{1}-2\sqrt{2}\alpha_{3})\pi+\frac{2505}{512}\alpha_{2}\pi+\frac{3195}{512}\alpha_{3}\pi),\\
\cdots\\
\end{aligned}
$$}

Following the lines of lemma 4.2, we obtain

\vskip 0.2 true cm
{\it \noindent{\bf Lemma 4.4.}~$M_{2}(h)$ has exactly 3 zeros near $h=0$.
}

\vspace{1em}
\noindent{\bf(2)~Homoclinic bifurcation of $l_{2}^{+}$}

For $n=3$, system (1.6) is equivalent to
$$
\begin{cases}
\dot{x}=y+2x^2y-2y^3,\\
\dot{y}=x-2x^3-2xy^2+\varepsilon(\bar{\alpha}_{0}y+\bar{\alpha}_{1}xy+\bar{\alpha}_{2}y^3+\bar{\alpha}_{3}x^2y),
\end{cases}\eqno(4.5)
$$
where
$$
\begin{tabular}{lll}
$\bar{\alpha}_{0}=-\frac{1}{2}\alpha_{0}+\frac{\sqrt{2}}{4}\alpha_{1}+\frac{3}{4}\alpha_{3},$&\qquad\qquad&$\bar{\alpha}_{1}=-\frac{1}{2}\alpha_{1}-\frac{\sqrt{2}}{2}\alpha_{3},$\\
\ &\qquad\qquad&\ \\
$\bar{\alpha}_{2}=-\frac{1}{2}\alpha_{2},$&\qquad\qquad&$\bar{\alpha}_{3}=-\frac{1}{2}\alpha_{3},$\\
\end{tabular}\eqno(4.6)
$$
with det$\frac{\partial(\bar{\alpha}_{0},\bar{\alpha}_{1},\bar{\alpha}_{2},\bar{\alpha}_{3})}{\partial(\alpha_{0},\alpha_{1},\alpha_{2},\alpha_{3})}=\frac{1}{16}$.
System (4.5)$|_{\varepsilon=0}$ has a double homoclinic orbit $L_{0}(=L_{10}\cup L_{20})$ given by the $H_{3}(x,y)=0$ with a hyperbolic saddle at the origin satisfying $H_{3}(0,0)=0$, where $L_{0}=\{(x,y)|H_{3}(x,y)=0\},$
$$L_{10}=\{(x,y)|H_{3}(x,y)=0,x>0\},\quad L_{20}=\{(x,y)|H_{3}(x,y)=0,x<0\}$$
with
$$H_{3}(x,y)=\frac{1}{2}y^2-\frac{1}{2}x^2+x^2y^2-\frac{1}{2}y^4+\frac{1}{2}x^4.$$
Near $L_{0}$ there are three families of periodic orbits $L(h)$, which are located outside $L_{0}$ for $h\in(0,\frac18)$, and $L_{1}(h)$ and $L_{2}(h)$, which are located inside $L_{0}$ for $h\in(-\frac18,0)$.
Then, correspondingly, $I(h)$ can be expressed as
$$I_{j}(h)=\oint_{L_{j}(h)}\bar{\alpha}_{0}y+\bar{\alpha}_{1}xy+\bar{\alpha}_{2}y^3+\bar{\alpha}_{3}x^2ydx,\ j=1,2$$
when $h\in(-\frac18,0)$, and
$$I_{3}(h)=\oint_{L(h)}\bar{\alpha}_{0}y+\bar{\alpha}_{2}y^3+\bar{\alpha}_{3}x^2ydx$$
when $h\in(0,\frac18)$.

The homoclinic loop $L_{j0},i=1,2$ can be expressed as
$$L_{j0}:
\begin{cases}
y^{+}(x)=\frac{1}{2}\sqrt{2+4x^2-2\sqrt{8x^4+1}},\ \ \ \ \ y>0,\\
y^{-}(x)=-\frac{1}{2}\sqrt{2+4x^2-2\sqrt{8x^4+1}},\ \ \ y<0.
\end{cases}
$$
Set $g(x,y,\bar{\delta})=\bar{\alpha}_{0}y+\bar{\alpha}_{1}xy+\bar{\alpha}_{2}y^3+\bar{\alpha}_{3}x^2y$ with $\bar{\delta}=(\bar{\alpha}_{0},\bar{\alpha}_{1},\bar{\alpha}_{2},\bar{\alpha}_{3})$,
by the expansions of $I(h)$ near the double homoclinic loop through $(0,0)$ with $H_{3}(0,0)=0$ in \cite{hanyangli}, we get the following lemma.
\begin{figure}[htbp]\centering
\centering
\includegraphics[scale=0.43]{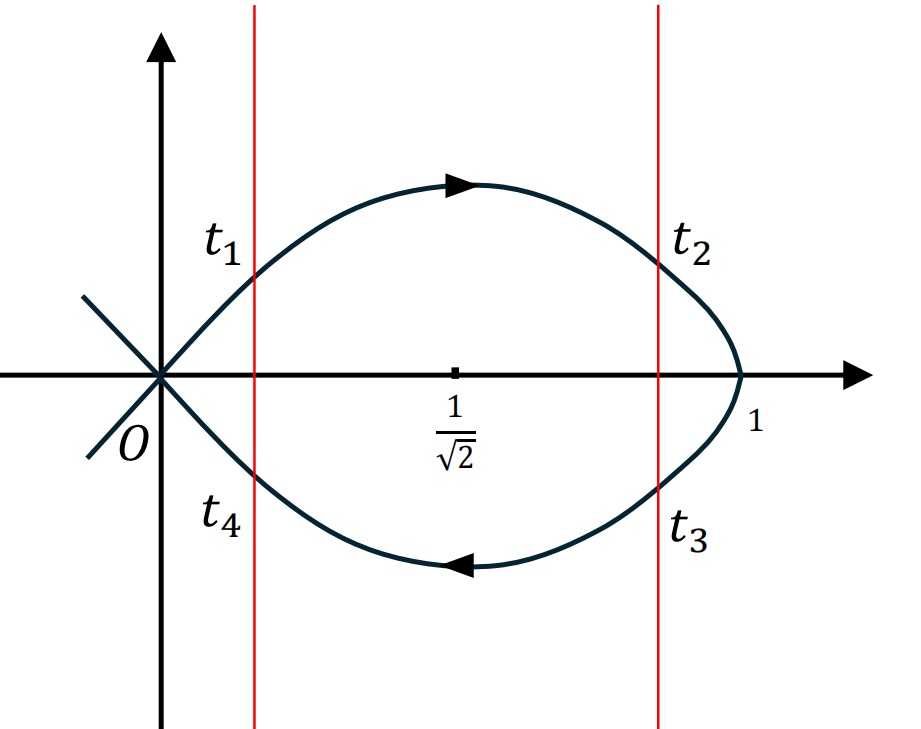}
\caption{\small About the calculation of $c_{2,1}$ and $c_{2,2}$: $x_{1}^{*}=x(t_{1})=x(t_{4})$, $x_{2}^{*}=x(t_{2})=x(t_{3})$}
\end{figure}

\vskip 0.2 true cm
{\it \noindent{\bf Lemma 4.5.}~Consider system (4.5), let $L_{j0}$ be expressed as $(x_{j}(t),y_{j}(t))$ for $t\in(-\infty,+\infty)$, $j=1,2$, we obtain

\noindent(i)~For $0<-h\ll1$,
$$I_{j}(h)=c_{0,j}+c_{1}hln|h|+c_{2,j}h+c_{3}h^2ln|h|+O(h^4),\ j=1,2,$$
for $0<h\ll1$,
$$I_{3}(h)=c_{0}+2c_{1}hln|h|+c_{2}h+2c_{3}h^2ln|h|+O(h^4),$$

\noindent(ii)~
$$
\begin{tabular}{lll}
$c_{0,1}=\bar{\alpha}_{0}A_{0}+\bar{\alpha}_{1}A_{1}+\bar{\alpha}_{2}A_{2}+\bar{\alpha}_{3}A_{3},$&\qquad\qquad&$c_{1}=\bar{\alpha}_{0},$\\
$\ $&\qquad\qquad&$\ $\\
$c_{0,2}=\bar{\alpha}_{0}A_{0}-\bar{\alpha}_{1}A_{1}+\bar{\alpha}_{2}A_{2}+\bar{\alpha}_{3}A_{3},$&\qquad\qquad&$c_{3}=\frac{1}{2}(\bar{\alpha}_{3}-3\bar{\alpha}_{2}),$\\
$\ $&\qquad\qquad&$\ $\\
$c_{0}=c_{0,1}+c_{0,2},$&\qquad\qquad&$c_{2}=c_{2,1}+c_{2,2},$\\
\end{tabular}$$
$$
c_{2,j}=\oint_{L_{j0}}g_{y}(x_{j}(t),y_{j}(t),\bar{\delta})-g_{y}(x_{j}(t),y_{j}(t),\bar{\delta})|_{x=y=0}dt+qc_{1},j=1,2,
$$
where $q\neq0$ is a constant, and
$$A_{0}=\oint_{L_{10}}ydx,A_{1}=\oint_{L_{10}}xydx,A_{2}=\oint_{L_{10}}y^3dx,A_{3}=\oint_{L_{10}}x^2ydx.$$
\noindent(iii)~ For $0<x_{1}^{*}<\frac{1}{\sqrt{2}}<x_{2}^{*}<1$, $y^{*}_{2}=y_{+}(x_{2}^{*})$(see Fig 3.), $c_{2,1}$ and $c_{2,2}$ can be expressed as
$$c_{2,1}|_{c_{1}=0}=J_{1}^{(1)}+J_{2}^{(1)}+J_{3}^{(1)}, \quad c_{2,2}|_{c_{1}=0}=J_{1}^{(2)}+J_{2}^{(2)}+J_{3}^{(2)},$$
where
$$
\begin{cases}
J_{1}^{(1)}=\int_{0}^{x_{1}^{*}}\frac{x}{y_{+}}\frac{2\bar{\alpha}_{1}+2\bar{\alpha}_{3}x}{(1+2x^2-2y^2_{+})}dx+\int_{0}^{x_{1}^{*}}\frac{6\bar{\alpha}_{2}y_{+}}{1+2x^2-2y^2_{+}}dx,\\
\ \\
J_{2}^{(1)}=2\int_{x^{*}_{1}}^{x_{2}^{*}}\frac{\bar{\alpha}_{1}x+3\bar{\alpha}_{2}y_{+}^2+\bar{\alpha}_{3}x^2}{y_{+}(1+2x^2-2y^2_{+})}dx,\\
\ \\
J_{3}^{(1)}=2\int_{0}^{-y_{2}^{*}}\frac{\bar{\alpha}_{1}x+3\bar{\alpha}_{2}y^2+\bar{\alpha}_{3}x^2}{x(1-2y^2-2x^2)}dy,
\end{cases}
$$
$$
\begin{cases}
J_{1}^{(2)}=\int_{0}^{x_{1}^{*}}\frac{-x}{y_{+}}\frac{2\bar{\alpha}_{1}-2\bar{\alpha}_{3}x}{(1+2x^2-2y^2_{+})}dx+\int_{0}^{x_{1}^{*}}\frac{6\bar{\alpha}_{2}y_{+}}{1+2x^2-2y^2_{+}}dx,\\
\ \\
J_{2}^{(2)}=2\int_{x^{*}_{1}}^{x_{2}^{*}}\frac{-\bar{\alpha}_{1}x+3\bar{\alpha}_{2}y_{+}^2+\bar{\alpha}_{3}x^2}{y_{+}(1+2x^2-2y^2_{+})}dx,\\
\ \\
J_{3}^{(2)}=2\int_{0}^{-y_{2}^{*}}\frac{-\bar{\alpha}_{1}x+3\bar{\alpha}_{2}y^2+\bar{\alpha}_{3}x^2}{x(1-2y^2-2x^2)}dy.
\end{cases}
$$
}
\noindent{\bf Proof}~According to the expansions in \cite{hanyangli}, we can get (i)--(ii)

Let us give the expression of $c_{2,j}|_{c_{1}=0}$, $j=1,2$. For $c_{2,1}|_{c_{1}=0}$, set some $-\infty<t_{1}<t_{2}<t^{*}<t_3<t_4<+\infty$ such that
$$0<x(t_1)<\frac{1}{\sqrt{2}}<x(t_{2})<x(t^{*})=1, x(t_{2})=x(t_{3}), x(t_1)=x(t_4), y_{+}(x(t^{*}))=0,$$
see Fig 3, let $x(t_1)=x^{*}_{1}$, $x(t_2)=x^{*}_{2}$, $y_{2}^{*}=y_{+}(x^{*}_{2})>0$, we have
$$
\begin{aligned}
c_{2,1}|_{c_{1}=0}=&\oint_{L_{10}}g_{y}(x_{1}(t),y_{1}(t),\bar{\delta})-g_{y}(x_{1}(t),y_{1}(t),\bar{\delta})|_{x=y=0}dt\\
=&\oint_{L_{10}}\omega(x_{1}(t),y_{1}(t))dt\\
=&J_{1}^{(1)}+J_{2}^{(1)}+J_{3}^{(1)},\\
\end{aligned}
$$
where $\omega(x_{1}(t),y_{1}(t))=\bar{\alpha}_{1}x_{1}(t)+3\bar{\alpha}_{2}y_{1}^2(t)+\bar{\alpha}_{3}x_{1}^2(t)$, and
$$
\begin{aligned}
J_{1}^{(1)}&=\left(\int_{-\infty}^{t_{1}}+\int_{t_{4}}^{+\infty}\right)\omega(x(t),y(t))dt,\\
J_{2}^{(1)}&=\left(\int_{t_{1}}^{t_{2}}+\int_{t_{3}}^{t_{4}}\right)\omega(x(t),y(t))dt,\\
J_{3}^{(1)}&=\int_{t_{2}}^{t_{3}}\omega(x(t),y(t))dt.\\
\end{aligned}
$$
Noting that $H_{3}(x,-y)=H_{3}(x,y)$, $\omega(x,-y)=\omega(x,y)$, we get
$$
\begin{aligned}
J_{1}^{(1)}&=\int_{0}^{x^{*}_{1}}\frac{\omega(x,y_{+})}{H_{3y}(x,y_{+})}dx+\int^{0}_{x^{*}_{1}}\frac{\omega(x,y_{-})}{H_{3y}(x,y_{-})}dx=2\int_{0}^{x_{1}^{*}}\frac{\bar{\alpha}_{1}x+3\bar{\alpha}_{2}y_{+}^2+\bar{\alpha}_{3}x^2}{y_{+}(1+2x^2-2y^2_{+})}dx\\
&=\int_{0}^{x_{1}^{*}}\frac{x}{y_{+}}\frac{2\bar{\alpha}_{1}+2\bar{\alpha}_{3}x}{(1+2x^2-2y^2_{+})}dx+\int_{0}^{x_{1}^{*}}\frac{6\bar{\alpha}_{2}y_{+}}{1+2x^2-2y^2_{+}}dx,\\
J_{2}^{(1)}&=\int_{x^{*}_{1}}^{x^{*}_{2}}\frac{\omega(x,y_{+})}{H_{3y}(x,y_{+})}dx+\int^{x^{*}_{2}}_{x^{*}_{1}}\frac{\omega(x,y_{-})}{H_{3y}(x,y_{-})}dx=2\int_{x^{*}_{2}}^{x_{1}^{*}}\frac{\bar{\alpha}_{1}x+3\bar{\alpha}_{2}y_{+}^2+\bar{\alpha}_{3}x^2}{y_{+}(1+2x^2-2y^2_{+})}dx.
\end{aligned}
$$
Noting that $x(-y)=x(y)$ for $x\in(x_{2}^{*},1)$, where
$$x(y)=\frac{1}{2}\sqrt{2-4y^2+2\sqrt{8y^4-8y^2+1}},\ y\in(-y^{*}_{2},y^{*}_{2}),$$
we have
$$J_{3}^{(1)}=\left(\int_{t_{2}}^{t^{*}}+\int_{t^{*}}^{t_{3}}\right)\omega(x(t),y(t))dt=\left(\int_{y_{2}^{*}}^{0}+\int_{0}^{-y_{2}^{*}}\right)\frac{\omega}{-H_{3x}}(x(y),y)dy.$$
Let $u=-y$, we have
$$\int_{y_{2}^{*}}^{0}\frac{\omega(x(y),y)}{-H_{3x}(x(y),y)}dy=\int_{0}^{-y_{2}^{*}}\frac{\omega(x(-u),-u)}{-H_{3x}(x(-u),-u)}du.$$
Since $H_{3}(x,-u)=H_{3}(x,u)$, thus
$$J_{3}^{(1)}=2\int_{0}^{-y_{2}^{*}}\frac{\omega(x(y),y)}{-H_{3x}(x(y),y)}dy=2\int_{0}^{-y_{2}^{*}}\frac{\bar{\alpha}_{1}x+3\bar{\alpha}_{2}y^2+\bar{\alpha}_{3}x^2}{x(1-2y^2-2x^2)}dy.$$
Then we can get the formula of $c_{2,1}|_{c_{1}=0}$, $c_{2,2}|_{c_{1}=0}$ can be obtained similarly. This ends the proof. $\diamondsuit$

Hence, for homoclinic bifurcation of system (4.5), we have the following lemma

{\it \noindent{\bf Lemma 4.6.}~For $I_{j}(h)$, $j=1,2,3$, we obtain

\noindent(i)~If $\bar{\alpha}_{1}\neq0$, $I_{j}(h)$ has at most $3$ zeros for $0<-h\ll1$, $j=1,2$.

\noindent(ii)~$I_{3}(h)$ has at most $2$ zeros for $0<h\ll1$.}

\noindent{\bf Proof}~
Set $x_{1}^{*}=0.05$, $x_{2}^{*}=0.9$, and $\frac{x}{y_{+}(x)}\approx1$ for $x\in(0,0.05)$, we obtain by Matlab that
$$
A_{0}\approx0.5301166457,A_{1}\approx0.2939666274,A_{2}\approx0.0543804979,A_{3}\approx0.1912645804,
$$
and
$$
\begin{aligned}
&c_{2,1}|_{c_{1}=0}\approx\bar{\alpha}_{1}A_{4}+\bar{\alpha}_{2}A_{5}+\bar{\alpha}_{3}A_{6},\\
&c_{2,2}|_{c_{1}=0}\approx-\bar{\alpha}_{1}A_{4}+\bar{\alpha}_{2}A_{5}+\bar{\alpha}_{3}A_{6},\\
\end{aligned}
$$
with
$$A_{4}\approx4.4879224539,A_{5}\approx1.1424442577,A_{6}\approx2.8725107531.$$

Set
$$
\begin{cases}
\mu_{1}=\bar{\alpha}_{1}A_{1}+\bar{\alpha}_{2}A_{2}+\bar{\alpha}_{3}A_{3},\\
\mu_{2}=\bar{\alpha}_{0},\\
\mu_{3}=q\bar{\alpha}_{0}+\bar{\alpha}_{1}A_{4}+\bar{\alpha}_{2}A_{5}+\bar{\alpha}_{3}A_{6},\\
\end{cases}
$$
then we get
$$
\begin{cases}
\bar{\alpha}_{1}\approx-0.6787660039\bar{\alpha}_3+12.44695659\mu_{1}+(0.5924767814q-6.598338879)\mu_{2}-0.5924767814\mu_{3},\\
\bar{\alpha}_{2}\approx0.1520760752\bar{\alpha}_3-48.89601885\mu_{1}+(-3.202773199q+25.92059350)\mu_{2}+3.202773199\mu_{3},\\
\end{cases}
$$
and
$$I_{1}(h)=\mu_{1}+\mu_{2}hln|h|+\mu_{3}h+c_{3}^{*}h^2ln|h|+O(h^4)$$
with
$$c_{3}^{*}\approx0.2718858872\bar{\alpha}_3+73.34402828\mu_{1}+(4.804159798q-38.88089025)\mu_{2}-4.804159798\mu_{3}.$$
Denote $\mu=(\mu_{1},\mu_{2},\mu_{3})$. Since
$$\lim_{\mu\to 0}c_{3}^{*}\approx0.2718858872\bar{\alpha}_3,$$
for $\bar{\alpha}_3=\bar{\alpha}_{3}^{*}\neq0$, there exist a neighborhood $U$ of $\mu=0$ such that $c_{3}^{*}\neq0$ for $\mu\in U$. Set some $\mu^{*}=(\mu_{1}^{*},\mu_{2}^{*},\mu_{3}^{*})\in U$ such that $0\ll|\mu_{1}^{*}|\ll|\mu_{2}^{*}|\ll|\mu_{3}^{*}|$, then we get the conclusion (i) for $j=1$. The conclusion of $I_{2}(h)$ and $I_{3}(h)$ can be proved similarly. This ends the proof. $\diamondsuit$

\vspace{1em}
\noindent{\bf(3)~Coexistence of limit cycles for $l_{2}^{+}$}

\vspace{1em}
{\it \noindent{\bf Lemma 4.7.}~Let $(N_{M_1},N_{M_2},N_{I_{1}},N_{I_{2}},N_{I_{3}})$ denote the zero distribution of $M_{1}(h)$, $M_{2}(h)$, $I_{1}(h)$, $I_{2}(h)$ and $I_{3}(h)$, we can prove there exist some $(\alpha_{0},\alpha_{1},\alpha_{2},\alpha_{3})$ such that $M_{i}(h)$ and $I_{j}(h)$ has totally at least $3$ zeros with the following distributions.
$$
\begin{aligned}
&(3,0,0,0,0),(0,3,0,0,0),(0,0,3,0,0),(0,0,0,3,0),(1,2,0,0,0),(2,1,0,0,0),\\
&(2,0,1,0,0),(2,0,0,1,0),(0,2,0,1,0),(0,2,1,0,0),(1,1,1,0,0),(1,1,0,1,0),\\
&(1,0,1,1,0),(0,1,1,1,0),(1,0,0,0,2),(0,1,0,0,2),(0,0,1,1,2),(0,0,2,2,1).\\
\end{aligned}$$
}

{\bf\noindent Proof}~Let us consider $(0,0,2,2,1)$ first. Following the coefficients formula in lemma 4.4 and (4.6), we set
$$
\begin{cases}
\mu_{1}=(-\frac{1}{2}\alpha_{0}+\frac{\sqrt{2}}{4}\alpha_{1}+\frac{3}{4}\alpha_{3})A_{0}-(\frac{1}{2}\alpha_{1}+\frac{\sqrt{2}}{2}\alpha_{3})A_{1}-\frac{1}{2}\alpha_{2}A_{2}-\frac{1}{2}\alpha_{3}A_{3},\\
\mu_{2}=(-\frac{1}{2}\alpha_{0}+\frac{\sqrt{2}}{4}\alpha_{1}+\frac{3}{4}\alpha_{3})A_{0}+(\frac{1}{2}\alpha_{1}+\frac{\sqrt{2}}{2}\alpha_{3})A_{1}-\frac{1}{2}\alpha_{2}A_{2}-\frac{1}{2}\alpha_{3}A_{3},\\
\mu_{3}=-\frac{1}{2}\alpha_{0}+\frac{\sqrt{2}}{4}\alpha_{1}+\frac{3}{4}\alpha_{3}.\\
\end{cases}
$$
Then we get
$$
\begin{cases}
\alpha_{0}=0.5\alpha_3-2\mu_3-1.700873342-2.405398148\mu_1+2.405398148\mu_2,\\
\alpha_{1}=-1.414213562\alpha_3-3.401746684\mu_1+3.401746684\mu_2,\\
\alpha_{2}=-3.517153903\alpha_3-18.38894528\mu_1-18.38894528\mu_2+19.49657197\mu_3,\\
\end{cases}
$$
and
$$
\begin{aligned}
I_{1}(h)&=\mu_1+\mu_3hln|h|+c_{2,1}^{*}h+O(h^2ln|h|),\\
I_{2}(h)&=\mu_2+\mu_3hln|h|+c_{2,2}^{*}h+O(h^2ln|h|),\\
I_{3}(h)&=(\mu_{1}+\mu_{2})+2\mu_{3}hln|h|+c_2^{*}h+O(h^2ln|h|)
\end{aligned}
$$
with
$$
\begin{aligned}
c_{2,1}^{*}&\approx0.572820764\alpha_3+18.13756013\mu_1+2.870784807\mu_2+(q-11.13687335)\mu_3,\\
c_{2,2}^{*}&\approx0.572820764\alpha_3+2.870784807\mu_1+18.13756013\mu_2+(q-11.13687335)\mu_3,\\
c_{2}^{*}&\approx1.145641527\alpha_3+21.00834495\mu_1+21.00834495\mu_2-22.27374669\mu_3.
\end{aligned}
$$
Meanwhile, for $M_{1}(h)$ and $M_{2}(h)$, we have
$$
\begin{aligned}
b_{1}(\delta)&=\alpha_{0}\pi\approx1.570796327\alpha_3-7.556781151\mu_1+7.556781151\mu_2-6.283185310\mu_3,\\
d_{1}(\delta)&=(\alpha_{0}-2\sqrt{2}\alpha_{1}+2\alpha_{3})\pi\approx20.42035225\alpha_3+22.67034345\mu_1-22.67034345\mu_2-6.283185308\mu_3,
\end{aligned}
$$
Denote $\mu=(\mu_{1},\mu_{2},\mu_{3})$. Since
$$
\begin{aligned}
&\lim_{\mu\to 0}c_{2,1}\approx0.572820764\alpha_3,\quad\lim_{\mu\to 0}b_{1}(\delta)\approx1.570796327\alpha_3\\
&\lim_{\mu\to 0}c_{2,2}\approx0.572820764\alpha_3,\quad\lim_{\mu\to 0}d_{1}(\delta)\approx20.42035225\alpha_3\\
&\lim_{\mu\to 0}c_{2}\approx1.145641527\alpha_3,\\
\end{aligned}
$$
for $\alpha_3=\alpha_{3}^{**}\neq0$, there exist a neighborhood $U$ of $\mu=0$ such that $c_{2,1}$, $c_{2,2}$, $c_{2}$, $b_{1}(\delta)$ and $d_{1}(\delta)$ are not vanished. Then system (1.6) has no limit cycles by Hopf bifurcation with $\mu\in U$. If $\alpha_{1}=-\sqrt{2}\alpha_{3}$, we have $\mu_{1}=\mu_{2}$. Noting that $0<-h\ll1$ for $I_{1}(h)$, $I_{2}(h)$, and $0<h\ll1$ for $I_{3}(h)$, set some $\mu^{*}=(\mu_{1}^{*},\mu_{1}^{*},\mu_{3}^{*})\in U$ such that $0\ll|\mu_{1}^{*}|\ll|\mu_{3}^{*}|$, then we can obtain the distribution $(0,0,2,2,1)$. Other distribution cases can be proved similarly. This ends the proof. $\diamondsuit$

\vspace{1em}
{\bf\noindent Proof of Theorem 1.2}~Following the line of lemma 4.2, lemma 4.4, lemma 4.6 and lemma 4.7, we can prove Theorem 1.2 is ture for system (1.6). This ends the proof. $\diamondsuit$

\vspace{1em}
\noindent{\bf Acknowledgments.}
The authors would like to express their sincere appreciation to the reviewer for his/her helpful comments which helped with improving the quality of this work.
\vskip 0.5 true cm

\noindent{\bf Data Availability Statement.} Data sharing not applicable to this article as no datasets were generated or analysed during the current study.
\vskip 0.5 true cm

\noindent{\bf Conflict of Interest.} The authors declare that they have no conflict of interest.


\begin{thebibliography}{99}
\bibitem{Arnold}
V.I. Arnold, Ten problems in: theory of singularities and its applications. Adv. Sov. Math. 1 (1990) 1-8.
\bibitem{chenyu}
Y. Chen, J. Yu, The study on cyclicity of a class of cubic systems. Discrete $\&$ Continuous Dynamical Systems-Series B. 27 (2022) 6233-6256.
\bibitem{Chow}
S. Chow, C. Li, D. Wang,  Normal forms and bifurcation of planar vector fields. Cambridge University Press, (1994).
\bibitem{collin}
C.B. Collins, Static stars: Some mathematical curiosities. Journal of Mathematical Physics 18 (1977) 1374-1377.
\bibitem{Dumortier1}
F. Dumortier, C. Li, Perturbations from an elliptic Hamiltonian of degree four: I. Saddle loop and two saddle cycle. Journal of Differential Equations 176 (2001) 114-157.
\bibitem{Dumortier2}
F. Dumortier, C. Li, Perturbations from an elliptic Hamiltonian of degree four: II. Cuspidal loop. Journal of Differential Equations 175 (2001) 209-243.
\bibitem{Dumortier3}
F. Dumortier, C. Li, Perturbation from an elliptic Hamiltonian of degree four--III. global centre. Journal of differential equations 188 (2003) 473-511.
\bibitem{Dumortier4}
F. Dumortier, C. Li, Perturbation from an elliptic Hamiltonian of degree four--IV. figure eight-loop. Journal of Differential Equations 188 (2003) 512-554.
\bibitem{Gavrilov}
L. Gavrilov, I.D. Iliev, Quadratic perturbations of quadratic codimension-four centers. Journal of mathematical analysis and applications, 357 (2009) 69-76.
\bibitem{hanyangli}
M. Han, J. Yang, J. Li, General study on limit cycle bifurcation near a double homoclinic loop. Journal of Differential Equations, 347 (2023) 1-23.
\bibitem{horo}
E. Horozov, I.D. Iliev, Linear estimate for the number of zeros of Abelian integrals with cubic Hamiltonians. Nonlinearity 11 (1998) 1521-1537.
\bibitem{Iliev}
I.D. Iliev, C. Li, J. Yu, On the cubic perturbations of the symmetric 8-loop Hamiltonian. Journal of Differential Equations 269 (2020) 3387-3413.
\bibitem{khovan}
A.G. Khovansky, Real analytic manifolds with finiteness properties and complex Abelian integrals. Functional Analysis and its Applications 18 (1984) 119-128.
\bibitem{lichengzhi2009}
C. Li, C. Liu, J. Yang, A cubic system with thirteen limit cycles. Journal of Differential Equations 246 (2009) 3609-3619.
\bibitem{zhangzhifen2002}
W. Li, Y. Zhao, C. Li, Z. Zhang, Abelian integrals for quadratic centers having almost all their orbits formed by quartics. Nonlinearity 15 (2002) 863-885.
\bibitem{cjliu8}
C. Liu, Estimate of the number of zeros of Abelian integrals for an elliptic Hamiltonian with figure--of--eight loop,
Nonlinearity 16 (2003) 1151-1163.
\bibitem{Liuxia}
X. Liu, M. Han, Bifurcation of limit cycles by perturbing piecewise Hamiltonian systems. International Journal of Bifurcation and Chaos, 20 (2010) 1379-1390.
\bibitem{3body}
K. Meyer, G. Hall, D. Offin, Introduction to Hamiltonian dynamical systems and the n-body problem. Springer New York (2009).
\bibitem{Petrov1}
G.S. Petrov, Complex zeros of an elliptic integral. Functional Analysis and Its Applications 21 (1987) 247-248.
\bibitem{Petrov2}
G.S. Petrov, Complex zeros of an elliptic integral. Functional Analysis and its Applications 23 (1989) 160-161.
\bibitem{Petrov3}
G.S. Petrov, Non-oscillation of elliptic integrals. Functional Analysis and Its Applications 24 (1990) 205-210.
\bibitem{Petrov4}
G.S. Petrov, On the nonoscillation of elliptic integrals. Functional Analysis and Its Applications 31 (1997) 262-265.
\bibitem{puu}
T. Puu, Attractors, bifurcations, $\&$ chaos: Nonlinear phenomena in economics. Springer Science $\&$ Business Media (2013).
\bibitem{varchen}
A.N. Varchenko, Estimate of the number of zeros of an Abelian integral depending on a parameter and limit cycles. Functional Analysis and its Applications 18 (1984) 98-108.
\bibitem{Wujuanjuan}
J. Wu, Y. Zhang, C. Li, On the number of zeros of Abelian integrals for a kind of quartic Hamiltonians. Applied Mathematics and Computation 228 (2014) 329-335.
\bibitem{yang2019ds}
J. Yang, On the number of zeros of Abelian integral for a class of cubic Hamiltonian systems. Dynamical Systems. 34 (2019) 561-583.
\bibitem{yangbf}
J. Yang, S. Sui, L. Zhao, On the Number of Zeros of Abelian Integral for a Class of Cubic Hamilton Systems with the Phase Portrait ``Butterfly''. Qualitative theory of dynamical systems. 18 (2019) 947-967.
\bibitem{yangjihua2017}
J. Yang, L. Zhao, The cyclicity of period annuli for a class of cubic Hamiltonian systems with nilpotent singular points. Journal of Differential Equations 263 (2017) 5554-5581.
\bibitem{zhaoyulin1999}
Y. Zhao, Z. Zhang, Linear Estimate of the number of zeros of Abelian integrals for a kind of quartic Hamiltonians. Journal of Differential Equations 155 (1999) 73-88.
\bibitem{zhouli1}
X. Zhou, C. Li, Estimate of the number of zeros of Abelian integrals for a kind of quartic Hamiltonians with two centers. Applied mathematics and computation 204 (2008) 202-209.
\bibitem{zhouli2}
X. Zhou, C. Li, On the algebraic structure of Abelian integrals for a kind of perturbed cubic Hamiltonian systems. Journal of mathematical analysis and applications 359 (2009) 209-215.
\bibitem{zhao2021}
J. Zhou, L. Zhao, J. Wang, Cyclicity of a class of Hamiltonian systems under perturbations of piecewise smooth polynomials. International Journal of Bifurcation and Chaos 31 (2021) 2150199.
\end{thebibliography}
\end{document}